\documentclass[11pt, reqno]{amsart}

\usepackage{amssymb}
\usepackage{hyperref} 
\usepackage{mathrsfs,bbold}
\usepackage{amsmath}
\usepackage{amsthm}
\usepackage{enumerate}
\usepackage{dsfont}
\usepackage{color}
\usepackage[a4paper]{geometry}
\usepackage[utf8]{inputenc}
\usepackage{stmaryrd}
\usepackage{titlesec}
\usepackage{mathabx}
\usepackage{appendix}
\usepackage{todonotes, fancyhdr}
\usepackage{chngcntr}
\usepackage{cancel}

\usepackage{setspace}
\renewcommand{\baselinestretch}{0.99}

\usepackage[all]{xy}
\numberwithin{subsection}{section}
\numberwithin{subsubsection}{subsection}
\numberwithin{equation}{section} 

\newenvironment{Dem}[1][\unskip]{%
    \begin{list}{\hspace{0.5cm}{\sf \textbf{Proof #1 --}}}{%
        \setlength{\topsep}{0pt}%
        \setlength{\leftmargin}{0pt}%
        \setlength{\rightmargin}{0pt}%
        \setlength{\listparindent}{0pt}%
        \setlength{\itemindent}{0pt}%
        \setlength{\parsep}{0pt}%
        \addtolength{\leftmargin}{20pt}%
        \addtolength{\rightmargin}{0pt}%
    } \item }{\hfill $\rhd$\end{list}\smallskip}

\renewcommand\thesection       {\arabic{section}}
\renewcommand\thesubsection    {\thesection{\boldmath $.$}\arabic{subsection}}
\renewcommand\thesubsubsection    {\thesection{\boldmath $.$}\arabic{subsection}{\boldmath $.$}\arabic{subsubsection}}

\titleformat{\section}[block]
{\filcenter\normalfont\sffamily\bfseries\Large}
{{\hspace{-0.7cm}}\thesection \hspace{0.2em} --\vspace{0.3cm}}{0.5em}{}

\titleformat{\subsection}[block]
{\filcenter\normalfont\sffamily\bfseries\large}  						  
{\hspace{-0.7cm}\thesubsection \hspace{0.5em} \vspace{0.3cm}}{.5em}{}  
\titlespacing{\subsection}{-0pc}{1.5ex plus .1ex minus .2ex}{0pc}

\titleformat{\subsubsection}[block]
{\filcenter\normalfont\sffamily\bfseries}					  
{\hspace{-0.7cm}\thesubsubsection \hspace{0.5em} \vspace{0.3cm}}{.5em}{}  
\titlespacing{\subsection}{-0pc}{1.5ex plus .1ex minus .2ex}{0pc}

\newtheoremstyle{mystyle}
{3pt}               
{3pt}               
{\it }                      
{}                      
{\sffamily\bfseries}             
{}                      
{0.5em}                 
{#1 #2{\Large$.$}  }

\theoremstyle{mystyle}

\newtheorem{thm}{Theorem}
\newtheorem*{thm*}{Theorem}

\newtheorem{cor}[thm]{\hspace{-0.15cm}  {Corollary} }
\newtheorem{lem}[thm]{\hspace{-0.2cm}  {Lemma} }
\newtheorem{prop}[thm]{\hspace{-0.13cm} {Proposition}}
\newtheorem{defn}[thm]{ \hspace{-0.25cm} {Definition}}

\newtheoremstyle{mystyle2}
{3pt}               
{3pt}               
{\it }                      
{}                      
{\sffamily\bfseries}             
{}                      
{0.5em}                 
{\llap{#2 }#1{\hspace{0.2cm}--}}

\theoremstyle{mystyle2}

\newtheorem*{definition*}{Definition}
\newtheorem*{theorem*}{Theorem}
\newtheorem*{Remark*}{Remark}
\newtheorem*{lem*} {Lemma}
\newtheorem*{defn*} {Definition}
\newtheorem*{prop*} {Proposition}
\newtheorem*{cor*} {Corollary}

\newcommand{\Nat}{\mathbb{N}}
\newcommand{\Int}{\mathbb{Z}}
\newcommand{\Euc}{\mathbb{R}}
\newcommand{\Hoel}{\mathcal{C}}
\newcommand{\Schw}{\mathcal{S}}

\newcommand{\imunit}{\sqrt{-1}}

\newcommand{\innpro}[2]{{\langle#1,#2\rangle}}

\newcommand{\iden}{\text{\rm Id}}

\newcommand{\supp}[1]{\mathop{\text{\rm supp}}(#1)}

\newcommand{\Fouri}{\mathscr{F}}

\newcommand{\unit}{\mathbf{1}}

\newcommand{\admi}{\mathcal{J}}
\newcommand{\remind}{\mathcal{N}}

\newcommand{\copro}{\Delta^+}
\newcommand{\comul}{\Delta}

\newcommand{\Brac}[1]{[\![#1]\!]}

\newcommand{\trino}[1]{|\!|\!|#1|\!|\!|}
\newcommand{\brarb}[1]{[\hspace{-0.5mm}]#1[\hspace{-0.5mm}]}

\newcommand{\bsf}{\boldsymbol{f}}
\newcommand{\bsh}{\boldsymbol{h}}
\newcommand{\bsX}{\boldsymbol{X}}

\newcommand{\ssk}{\smallskip}

\renewcommand{\epsilon}{\varepsilon}

\newcommand\bbN{\mathbb{N}}
\newcommand\bbR{\mathbb{R}}

\newcommand\bbZ{\mathbb{Z}}

\newcommand{\mcA}{\mathcal{A}}
\newcommand{\mcB}{\mathcal{B}} 
\newcommand{\mcC}{\mathcal{C}} 
\newcommand{\mcD}{\mathcal{D}}

\newcommand{\mcI}{\mathcal{I}}
\newcommand{\mcJ}{\mathcal{J}}

\newcommand{\mcN}{\mathcal{N}}

\newcommand\mcS{\mathcal{S}}

\newcommand{\bfD}{\mathbf{D}}

\newcommand{\bfI}{\mathbf{I}}

\newcommand{\bfP}{\mathbf{P}}
\newcommand{\bfQ}{\mathbf{Q}}

\newcommand{\bfX}{\mathbf{X}}

\begin{document}

\begin{center}
{\LARGE\sffamily{Paracontrolled calculus and regularity structures (1)  \vspace{0.5cm}}}
\end{center}

\begin{center}
{\sf I. BAILLEUL} \& {\sf M. HOSHINO}
\end{center}

\vspace{1cm}

\begin{center}
\begin{minipage}{0.8\textwidth}
\renewcommand\baselinestretch{0.7} \scriptsize \textbf{\textsf{\noindent Abstract.}} We start in this work the study of the relation between the theory of regularity structures and paracontrolled calculus. We give a paracontrolled representation of the reconstruction operator and provide a natural parametrization of the space of admissible models.
\end{minipage}
\end{center}

\vspace{0.6cm}

\section{Introduction}

Starting with his groundbreaking work \cite{Hai14}, M. Hairer has developed with his co-authors \cite{BHZ16, CH16, BCCH18} a theory of subcritical singular stochastic partial differential equations (PDEs) that provides now an automated blackbox for the basic understanding of a whole class of stochastic PDEs. Equations of this class all share the common feature of involving ill-defined products of distributions with functions or distributions. The methodology of regularity structures for the study of a given singular stochastic PDE takes its roots in T. Lyons' theory of rough paths, such as reshaped by M. Gubinelli \cite{GubinelliControlled, GubinelliBranched}. It requires first to identify a proper space of enhanced noises. The raw random noise that appears in the equation needs to be lifted into a random noise taking values in that enhanced space. This is typically a probabilistic task, mostly independent of the details of the dynamics under study, once the appropriate space of enhanced noises has been constructed from the equation. (That space happens to be equation-independent in the rough differential equation setting, while it is equation-dependent in a PDE setting.) The lifting task typically involves stochastic or Gaussian calculus in a rough paths setting; it involves the difficult implementation of a renormalisation procedure in the singular stochastic PDE setting. This step somehow takes care of the core problem: defining the product of two random distributions as a random variable rather than taking the product of two realizations of these random variables. \textit{These enhanced noises come under the form of a model in regularity structures}. This is a deterministic object, and the previous step takes care of constructing a random model. Having a model is somewhat equivalent to having a definition of the product of a number of otherwise possibly ill-defined quantities. A restricted class of space-time functions or distributions is then described in regularity structures theory under the form of a space-time indexed family of jets describing them locally around each space-time point. Given any choice of model, a consistency relation ensures that coherent jets describe indeed true space-time functions or distributions. This is the role of the reconstruction operator; \textit{coherent jets are modelled distributions}. It happens then that one can reformulate the formal ill-posed equation into the space of jets as a well-posed, model-dependent, fixed point equation in a well-chosen space of jets. For the random model built from a renormalisation procedure in \cite{CH16}, the space-time function/distribution associated with the solution of the fixed point equation on the jet space can be shown to be the limit in probability of solutions of a family of well-posed space-time stochastic PDEs driven by regularized noises, as the regularization parameter tends to $0$ -- this is the content of \cite{BCCH18}. The fact that some of the terms in these modified and regularized stochastic PDEs blow up as the regularization parameter goes to $0$ is a feature of the singular nature of the initial equation.

Let us emphasize that the multiplication problem is fundamentally dealt with on the ground of the following heuristic argument. If one can make sense of the product of a number of reference quantities, one can make sense of the product of quantities that look like the reference quantities. This is what motivates the introduction of jets on scene.

\smallskip

The choice of a jet space to describe a possible solution to a singular stochastic PDE is not the only possible. As a matter of fact, Gubinelli, Imkeller and Perkowski devised in \cite{GIP15} a Fourier-based approach to the study of  singular stochastic PDEs whose scope has been extended in \cite{BB15, BBF15, BB16}. The heuristic remains the same, but paraproducts are used as a mean of making sense of what it means to look like a reference distribution or function. This choice of representation makes the technical details of paracontrolled calculus rather different from their regularity structures counterparts, and paracontrolled calculus remains to be systematized. Despite that fact, it happens to be possible to make a close comparison between the two settings. We start that comparison in this work by providing an 'explicit' paracontrolled representation of the reconstruction operator. This is the operator that associates to a coherent jet a space-time distribution. All notions and notations in the statement are properly defined below.

\medskip

\begin{thm}
\label{ThmMain1}
Let a concrete regularity structure $\mathscr{T}=(T^+,T)$ be given, together with a model $\sf M=(g,\Pi)$ on it.   \vspace{0.1cm}

\begin{enumerate}
   \item[\emph{\sf (1)}] One can construct functions $\Brac{\cdot}^{\sf M} : T\mapsto \mcC^{\beta_0}(\bbR^d)$ and $\Brac{\cdot}^{\sf g} : T^+\mapsto \mcC^0(\bbR^d)$, such that \vspace{0.1cm}
   
   \begin{itemize}
		\item[--] $\Brac{\sigma}^{\sf M}\in \mcC^{\vert\sigma\vert}(\bbR^d)$, and $\Brac{\tau}^{\sf g}\in \mcC^{\vert\tau\vert}(\bbR^d)$, for every homogeneous $\sigma\in T$ and $\tau\in T^+$,   \vspace{0.1cm}
 
		\item[--] all $\Brac{\sigma}^{\sf M}$, and $\Brac{\tau}^{\sf g}$ are continuous function of the model $\sf (g,\Pi)$,   \vspace{0.1cm}
   \end{itemize}  
and the following holds true.   \vspace{0.15cm}
   
   \item[\emph{\sf (2)}] One can associate to any modelled distribution 
$$
\bsf=\sum_{\tau\in\mcB; |\tau|<\gamma} f^\tau \tau \in \mcD^\gamma(\mathscr{T}, {\sf g}),
$$ 
a distribution $\Brac{\bsf}^{\sf M} \in \mcC^\gamma(\bbR^d)$ such that one defines a reconstruction $\textsf{\textbf{R}}\bsf$ of $\bsf$ setting
\begin{align}\label{2: eq reconst to PD}
\textsf{\textbf{R}}\bsf  :=  \sum_{\tau\in\mcB; |\tau|<\gamma} {\sf P}_{f^\tau} \Brac{\tau}^{\sf M} + \Brac{\bsf}^{\sf M}.
\end{align}
Each coefficient $f^\tau$, also has a representation
\begin{align}\label{2: eq coefficients to PD}
f^\tau  =  \sum_{\tau<\mu; |\mu|<\gamma} {\sf P}_{f^\mu} \Brac{\mu/\tau}^{\sf g} + \Brac{f^\tau}^{\sf g},
\end{align}
for some $\Brac{f^\tau}^{\sf g}\in\mcC^{\gamma-|\tau|}(\bbR^d)$. Moreover, the map 
$$
\bsf \mapsto \Big(\Brac{\bsf}^{\sf M}, \big(\Brac{f^\tau}^{\sf g}\big)_{\tau\in\mcB} \Big)
$$
from $\mcD^\gamma(\mathscr{T}, {\sf g})$ to $\mcC^\gamma(\bbR^d)\times \prod_{\tau\in\mcB} \mcC^{\gamma-|\tau|}(\bbR^d)$, is continuous.
\end{enumerate}
\end{thm}

\smallskip

This is the content of Proposition \ref{PropDefnBracket} and Theorem \ref{ThmReconstructionPC}. Any regularity exponent $a\in\bbR$ is allowed in the above statement. The inductive definition of $\Brac{\cdot}^{\sf M}$, Proposition \ref{PropDefnBracket}, will make it clear that $\Brac{\sigma}^{\sf M}$ can be understood as the `part' of $\sf \Pi\sigma$ of regularity $\mcC^{\vert\sigma\vert}(\bbR^d)$. The quantity $\Brac{\tau}^{\sf g}$ has a similar meaning for the function ${\sf g}(\tau)$. Theorem \ref{ThmMain1} provides a much refined version of the paraproduct-based construction of the reconstruction operator from Gubinelli, Imkeller and Perkowski' seminal work \cite{GIP15}. Notice that this statement is not related to any problem about singular stochastic PDE. The treatment of such equations involves the additional ingredient of an abstract integration operator and the additional notion of admissible model. We provide an explicit paracontrolled-based parametrization of that set of models under some canonical structure assumptions on the regularity structure.

\smallskip

\begin{thm}
\label{ThmMain2}
Given any family of distributions $\big(\Brac{\tau}\in\mcC^{|\tau|}(\bbR^d)\big)_{\tau\in\mcB; |\tau|\leq 0}$, there exists a unique admissible model $\sf M=(g,\Pi)$ on $\mathscr{T}$ such that one has 
\begin{equation}
\label{EqStructurePiBracketIntro}
{\sf \Pi}\tau := \sum_{\sigma<\tau}{\sf P}_{{\sf g}(\tau/\sigma)}\Brac{\sigma} + \Brac{\tau},
\end{equation}
for all $\tau\in\mcB$ with $|\tau|\leq 0$.
\end{thm}

\smallskip

The fact that identity \eqref{EqStructurePiBracketIntro} holds true with $\Brac{\cdot}^{\sf M}$ in place of $\Brac{\cdot}$ for \textit{any} model $\sf M=(g,\Pi)$, is part of the proof of item {\sf (1)} of Theorem \ref{ThmMain1}. 

\medskip

We work throughout with the usual isotropic H\"older spaces. All the results presented here have direct analogues involving anisotropic H\"older spaces, such as required for applications to parabolic singular stochastic PDEs. The proofs of all results are strictly identical. We refrain from putting ourselves in that setting so as not to overload the reader with additional technical details and keep focused on the main novelty. The reader will find relevant technical details in the work \cite{MP18} of Martin and Perkowski.

\bigskip

No previous knowledge of regularity structures or paracontrolled calculus is needed in this work, that is mostly self-contained, with the exception of elementary facts on paraproducts recalled in Appendix \ref{AppendixParapdcts}. We have thus given at few places full proofs of statements that were first proved elsewhere. The material has been organized as follows. Section \ref{SectionBasicsRS} sets the scene of regularity structures under a convenient form for us: Concrete regularity structures, models and modelled distributions are introduced, together with a number of elementary identities and examples. Theorem \ref{ThmMain1} is proved in Section \ref{SubsectionExplicitReconstruction}, while Section \ref{SectionParametrization} takes care of Theorem \ref{ThmMain2}.

\bigskip

\noindent \textsf{\textbf{Notations\, -- }}\textit{We use exclusively the letters $\alpha, \beta, \gamma, \theta$ to denote real numbers, and use the letters  $\sigma, \tau, \mu$ to denote elements of $T$ or $T^+$. We agree to use the shorthand notation $\frak{s}^{(+)}$ to mean both the statement $\frak{s}$ and the statement $\frak{s}^+$.}

\bigskip

\section{Basics on regularity structures}
\label{SectionBasicsRS}

Regularity structures are the backbone of expansion devices for the local description of functions and distributions in $\bbR^d$. The usual notion of local description of a function, near a point $x\in\bbR^d$, involves Taylor expansion and amounts to comparing a function to a polynomial centered at $x$
\begin{equation}   \label{EqLocalDescription}
f(\cdot) \simeq \sum_k f^k(x)\,(\cdot-x)^k, \quad \textrm{near }x.
\end{equation}
The sum over $k$ is finite and the approximation quantified. One gets a local description of $f$ near another point $x'$ writing
\begin{equation}   \textcolor{red}{\label{EqDescriptionf}}   \begin{split}
f(\cdot) \simeq \sum_{\ell\leq k} f^k(x)\,{k\choose{\ell}} (\cdot-x')^\ell (x'-x)^{k-\ell} \simeq \sum_\ell \left(\sum_{k;\ell\leq k} f^k(x)\,{k\choose{\ell}} (x'-x)^{k-\ell}\right)(\cdot-x')^\ell.
\end{split}\end{equation}
A more general local description device involves an $\bbR^d$-indexed collection of functions or distributions $(\Pi_x\tau)(\cdot)$, with labels in a finite set $\mcB=\{\tau\}$. Consider the real vector space $T$ spanned freely by $\mcB$. Functions or distributions are locally described as
$$
f(\cdot) \simeq \sum_\tau f^\tau(x)(\Pi_x\tau)(\cdot), \quad\textrm{near each }x\in\bbR^d.
$$
This implicitely assumes that the coefficients $f^\tau(x)$ are function of $x$. One has $\{\tau\}=\{k\}$ and $(\Pi_xk)(\cdot) = (\cdot - x)^k$, in the polynomial setting. Like in the former setting, in a general local description device the reference objects 
\begin{equation}   \textcolor{red}{\label{EqTransitionMapsRelation}}
(\Pi_{x'}\tau)(\cdot) = \big(\Pi_x(\Gamma_{xx'}\tau)\big)(\cdot)
\end{equation}
at a different base point $x'$ are linear combinations of the $\Pi_x\sigma$, for a linear map 
$$
\Gamma_{xx'} : T\rightarrow T, 
$$
and one can switch back and forth between local descriptions at different points. The linear maps $\Gamma_{xx'}$ are thus invertible and one has a group action of an $\bbR^d\times\bbR^d$-indexed group on the local description structure $T$.   

\smallskip

Whereas one uses the same polynomial-type local description for the $f^k$ as for $f$ itself in the usual $\mcC^\alpha$ setting, there is no reason in a more general local description device to use the same reference objects for $f$ and for its local coefficients, especially if the $(\Pi_x\tau)(\cdot)$ are meant to describe distributions, among others, while it makes sense to use functions only as reference objects to describe the functions $f^\tau$. A simple setting consists in having all the $f^\tau$ locally described by a possibly different finite collection $\mcB^+=\{\mu\}$ of labels, in terms of reference functions ${\sf g}_{yx}(\mu)$, with
$$
f^\tau(y) \simeq \sum_{\mu\in\mcB} f^{\tau\mu}(x) {\sf g}_{yx}(\mu), \quad\textrm{near }x.
$$
One thus has both
\begin{equation}   \label{EqGammaAppears}
f(\cdot) \simeq \sum_{\tau\in\mcB} f^\tau(x)(\Pi_x\tau)(\cdot) \simeq \sum_{\tau\in\mcB,\,\mu\in\mcB^+} f^{\tau\mu}(y) {\sf g}_{yx}(\mu)(\Pi_x\tau)(\cdot)
\end{equation}
and 
$$
f(\cdot) \simeq \sum_{\sigma\in \mcB} f^\sigma(y)(\Pi_y\sigma)(\cdot).
$$
\textit{Consistency} dictates that the two expressions coincide, giving in particular the fact that the coefficients $f^{\tau\mu}(y)$ are linear combinations of the $f^\sigma(y)$. Re-indexing identity \eqref{EqGammaAppears} and using the notation $\sigma/\tau$ for the $\mu$ corresponding to $\tau\mu\simeq\sigma$, one then has
\begin{equation}   \label{EqFirstExpansion}
f(\cdot) \simeq \sum_{\sigma\in \mcB, \tau\in\mcB} f^\sigma(y)\,{\sf g}_{yx}(\sigma/\tau)(\Pi_x\tau)(\cdot).
\end{equation}
The transition map $\Gamma_{xy} : T\rightarrow T$, from \eqref{EqTransitionMapsRelation} is thus given in terms of the splitting map 
$$
\Delta : T\rightarrow T\otimes T^+, \quad \Delta\sigma = \sum_{\tau} \tau\otimes (\sigma/\tau)
$$
that appears in the above decomposition, with
$$
\Pi_y\sigma = \sum_{\tau\in\mcB}{\sf g}_{yx}(\sigma/\tau)\Pi_x\tau
$$
so 
$$
\Gamma_{xy}\sigma = \sum_{\tau\in\mcB}{\sf g}_{yx}(\sigma/\tau) \tau.
$$
If one further expands $f^\sigma(y)$ in \eqref{EqFirstExpansion} around another reference point $z$, one gets
\begin{equation} \label{EqDeltaDeltaPlus} \begin{split}
f(\cdot) &\simeq \sum_{\tau,\sigma,\nu\in\mcB} f^\nu(z)\,{\sf g}_{zy}(\nu/\sigma){\sf g}_{yx}(\sigma/\tau)(\Pi_x\tau)(\cdot)   \\
		      & \simeq \sum_{\nu\in\mcB} f^\nu(z)(\Pi_z\nu)(\cdot) \simeq \sum_{\tau,\nu\in\mcB} f^\nu(z)\,{\sf g}_{zx}(\nu/\tau)(\Pi_x\tau)(\cdot).
\end{split} \end{equation}
Here again, \textit{consistency} requires that the two expressions coincide, giving the identity
$$
\sum_{\sigma\in\mcB} {\sf g}_{zy}(\nu/\sigma){\sf g}_{yx}(\sigma/\tau) = {\sf g}_{zx}(\nu/\tau)
$$
in terms of another splitting map
$$
\Delta^+ : T^+\rightarrow T^+\otimes T^+
$$
satisfying by construction the itendity
$$
(\textrm{Id}\otimes \Delta^+)\Delta = (\Delta\otimes \textrm{Id})\Delta
$$
encoded in identity \eqref{EqDeltaDeltaPlus}. Developing $f^\nu(z)$ in \eqref{EqDeltaDeltaPlus} in terms of another reference point leads by consistency to the identity
$$
(\textrm{Id}\otimes \Delta^+)\Delta^+ = (\Delta^+\otimes \textrm{Id})\Delta^+.
$$
If we insist that the family of reference functions ${\sf g}_{yx}(\mu), \mu\in\mcB^+$, be sufficiently rich to describe locally an \textit{algebra} of functions, it is convenient to assume that the linear span $T^+$ of $\mcB^+$ has an algebra structure and the maps ${\sf g}_{yx}$ on $T^+$ are characters of the algebra -- multiplicative maps. Building on the example of the polynomials, it is also natural to assume that $T^+$ has a grading structure; an elementary fact from algebra then leads directly to the Hopf algebra structure that appears below in the definition of a concrete regularity structure.   

\medskip

We choose to record the essential features of this discussion in the definition of a concrete regularity structure given below; this is a special form of the more general notion of regularity structure from Hairer' seminal work \cite{Hai14}. The reader should keep in mind that the entire algebraic setting can be understood at a basic level from the above consistency requirements on a given local description device. We refer the reader to Sweedler's book \cite{Sweedler} for an accessible reference on Hopf algebras. Given two statements $\frak{s}$ and $\frak{s}^+$, recall the convention that we agree to write $\frak{s}^{(+)}$ to mean both the statement $\frak{s}$ and the statement $\frak{s}^+$.

\bigskip

\hfill \textcolor{gray}{\textsf{\textbf{Concrete regularity structures}}}   \vspace{0.2cm}

\begin{defn*}
A \textsf{\textbf{concrete regularity structure}} $\mathscr{T}=(T^+,T)$ is the pair of graded vector spaces   \vspace{-0.1cm}
$$
T^+ =: \bigoplus_{\alpha\in A^+} T_\alpha^+, \qquad T =: \bigoplus_{\beta\in A} T_\beta   \vspace{-0.1cm}
$$
such that the following holds.   \vspace{0.1cm}
\begin{itemize}
   \item The index set $A^{+}\subset \bbR_+$ contains the point $0$, and $A^++A^+ \subset A^+$; the index set $A\subset\bbR$ is bounded below, and both $A^+$and $A$ have no accumulation points in $\bbR$. Set 
$$
\beta_0 := \min A.
$$  
   
   \item The vector spaces $T_\alpha^+$ and $T_\beta$ are finite dimensional. \vspace{0.1cm}
   
   \item The set $T^+$ is an algebra with unit $\bf 1$, with a Hopf structure with coproduct 
   $$
   \Delta^+ : T^+\rightarrow T^+\otimes T^+,
   $$
   such that $\Delta^+{\bf 1} = {\bf 1}\otimes {\bf 1}$, and, for $\tau\in T_\alpha^+$,
   \begin{equation}   \label{EqDefnDeltaPlus}
   \Delta^+\tau \in \left\{\tau\otimes {\bf 1} + {\bf 1}\otimes \tau + \sum_{0<\beta<\alpha} T^+_\beta\otimes T^+_{\alpha-\beta}\right\},
   \end{equation}
   \item One has $T_0^+ =\langle {\bf 1}\rangle$, and for any $\alpha,\beta\in A^+$, one has $T_\alpha^+ T_\beta^+ \subset T_{\alpha+\beta}^+$.   \vspace{0.1cm}
   
   \item One has a splitting map 
   $$
   \Delta : T \rightarrow T\otimes T^+,
   $$ 
   of the form
   \begin{equation}
   \label{EqDelta}  
   \Delta \tau \in \left\{\tau\otimes {\bf 1} + \sum_{\beta<\alpha} T_\beta\otimes T^+_{\alpha-\beta}\right\}   
   \end{equation}
   for each $\tau\in T_\alpha$, with the right comodule property
   \begin{equation}
   \label{EqRightComodule}   
   (\Delta\otimes \textrm{\emph{Id}})\Delta = (\textrm{\emph{Id}}\otimes \Delta^+)\Delta.
   \end{equation}
\end{itemize}
Let $\mcB_\alpha^+$  and $\mcB_\beta$ be bases of $T_\alpha^+$ and $T_\beta$, respectively. We assume $\mcB_0^+=\{{\bf1}\}$. Set 
$$
\mcB^+ := \bigcup_{\alpha\in A^+} \mcB_\alpha^+, \quad \mcB := \bigcup_{\beta\in A} \mcB_\beta.
$$ 
An element $\tau$ of $T_\alpha^{(+)}$ is said to be \emph{homogeneous} and is assigned \textsf{\textbf{homogeneity}} $|\tau| := \alpha$. The homogeneity of a generic element $\tau\in T^{(+)}$ is defined as $|\tau| := \max\{\alpha\}$, such that $\tau$ has a non-null component in $T^{(+)}_\alpha$.
We sometimes denote by 
$$
\mathscr{T} := \big((T^+,\Delta^+), (T,\Delta)\big)
$$
a concrete regularity structure.
\end{defn*}

\medskip

Note that we do not assume any relation between the linear spaces $T_\alpha^+$ and $T_\beta$ at that stage. Note also that the parameter $\beta$ in \eqref{EqDelta} can be non-positive, unlike in \eqref{EqDefnDeltaPlus}. For an arbitrary element $h$ in $T$, set
\begin{equation*}
h = \sum_{\beta\leq |h|} h_\beta \in \bigoplus_{\beta\leq |h|} T_\beta.
\end{equation*}
We use a similar notation for elements of $T^+$. For $\gamma\in\bbR$, set
\vspace{-0.1cm}$$
T_{<\gamma} := \bigoplus_{\beta<\gamma} T_\beta, \qquad T_{<\gamma}^+ := \bigoplus_{\alpha<\gamma} T_\alpha^+.
$$
The homogeneous spaces $T_{\beta}$ and $T^{+}_\alpha$ being finite dimensional, all norms on them are equivalent; we use a generic notation $\|\cdot\|_\beta$ or $\|\cdot\|_\alpha$ for norms on these spaces. For simplicity, we write
\begin{align}
\label{EqNotationDirectSum}
\|h\|_\alpha:=\|h_\alpha\|_\alpha.
\end{align}
To have a picture in mind, think of $T$ and $T^+$ as sets of possibly labelled rooted trees, with $T^+$ consisting only of trees with positive tree homogeneities -- a homogeneity is assigned to each labelled tree. This notion of homogeneity induces the decomposition \eqref{EqNotationDirectSum} of $T$ into linear spaces spanned by trees with the same homogeneity; a similar decomposition holds for $T^+$. The coproduct $\Delta^+\tau$ is typically a sum over subtrees $\sigma$ of $\tau$ with the same root as $\tau$, and $\tau/\sigma$ is the quotient tree obtained from $\tau$ by identifying $\sigma$ with the root. One understands the splitting $\Delta\tau$ of an element $\tau\in T$ in similar terms. See e.g. Section 2 and Section 3 of \cite{BHZ16}.   

\medskip

\noindent \textbf{\textsf{Notation.}} \textit{Given $\sigma, \tau \in \mcB^{(+)}$, we use the notation $\sigma\leq^{(+)}\tau$ to mean that $\sigma$ appears as a left hand side of one of the tensor products in the sum defining $\Delta^{(+)}\tau$; we write $\tau/^{(+)}\sigma$ for the corresponding right hand side, so we have, for $\tau\in\mcB^{(+)}$
$$
\Delta^{(+)}\tau = \sum_{\sigma\in\mcB^{(+)}} \sigma\otimes (\tau/^{(+)}\sigma).
$$
Write $\sigma<^{(+)}\tau$ to mean further that $\sigma$ is different from $\tau$. The notations $\tau/^{(+)}\sigma$ and $\sigma<^{(+)}\tau$ are only used for $\tau$ and $\sigma$ in $\mcB^{(+)}$.}

\medskip

Decomposing $\Delta\tau$ in the basis $\mcB\otimes\mcB^+$ of $T\otimes T^+$ as
$$
\Delta\tau =: \sum_{\sigma\in\mcB, \theta\in\mcB^+} (\Delta\tau)^{\sigma\theta} \, \sigma\otimes\theta,
$$
one has
\vspace{-0.1cm}$$
\tau/\sigma = \sum_{\theta\in\mcB^+} (\Delta\tau)^{\sigma\theta} \, \theta.
$$
We have a similar expression for $\tau/^+\sigma$; for $\sigma, \tau \in \mcB^+$,
\begin{equation}
\label{EqFormulaTauOnPlusSigma}
\tau/^+\sigma = \sum_{\theta\in\mcB^+} (\Delta^+\tau)^{\sigma\theta} \, \theta.
\end{equation}
With these notations, the right comodule property \eqref{EqRightComodule} writes for all $\tau\in\mcB$
\begin{equation}
\label{EqRightComoduleComponents}
\sum_{\sigma\in\mcB} (\Delta\tau)^{\sigma c}\,(\Delta\sigma)^{ab} = \sum_{\theta\in\mcB^+} (\Delta\tau)^{a\theta}\,(\Delta^+\theta)^{bc}
\end{equation}
for all $a\in\mcB$ and $b,c\in\mcB^+$. The identity from Lemma \ref{LemmaCoprodQuotient} is a direct consequence of the \textit{co-associativity property} 
$$
(\Delta^+\otimes \textrm{Id})\Delta^+ = ( \textrm{Id}\otimes \Delta^+)\Delta^+, 
$$
of the coproduct $\Delta^+$, and the right comodule identity \eqref{EqRightComodule}.   \vspace{0.1cm}

\begin{lem}
\label{LemmaCoprodQuotient}
For $\sigma<^+\tau$ in $\mcB^+$, we have
\begin{equation}
\label{EqCoprodQuotient+}
\begin{split}
\Delta^+(\tau/^+\sigma) &= \sum_{\sigma \leq^+ \eta \leq^+ \tau}(\eta/^+\sigma)\otimes(\tau/^+\eta)   \\
                                     &= (\tau/^+\sigma)\otimes\unit + \unit\otimes(\tau/^+\sigma) + \sum_{\sigma <^+ \eta <^+ \tau}(\eta/^+\sigma)\otimes(\tau/^+\eta).
\end{split}
\end{equation}
For $\sigma<\tau$ in $\mcB$, we have
\begin{equation}
\label{EqCoprodQuotient}
\begin{split}
\Delta^+(\tau/\sigma) &= \sum_{\sigma \leq \eta \leq \tau}(\eta/\sigma)\otimes(\tau/\eta).
\end{split}
\end{equation}
\end{lem}

\smallskip

A character $g$ on the algebra $T^+$ is a linear map $g:T^+\to\bbR$ such that $g(\tau_1\tau_2)=g(\tau_1)g(\tau_2)$ for any $\tau_1,\tau_2\in T^+$. The antipode $\mcA$ of the Hopf algebra structure turns the set of characters of the algebra $T^+$  into a group $G^+$ for the \textit{convolution law} $*$ defined by 
$$
(g_1*g_2)\tau := (g_1\otimes g_2)\Delta^+\tau, \quad \tau\in T^+.
$$
The identity of the group is the counit ${\bf 1}'$, the dual basis vector of the unit $\unit$, and the inverse $g^{-1} = g\circ\mcA$. One associates to a character $g$ of $T^+$ the map 
$$
\widehat g := (\textrm{Id}\otimes g)\Delta : T\mapsto T,
$$ 
from $T$ to itself. We have 
$$
\widehat{g_1*g_2} = \widehat{g_1} \circ \widehat{g_2}
$$
for any $g_1,g_2\in G^+$, as a consequence of the comodule property \eqref{EqRightComodule}. Also, for any $\tau\in T$,
$$
\Big(\widehat{g}(\tau) - \tau\Big) \in T_{<|\tau|},
$$
as a consequence of the structural identity \eqref{EqDelta}. Remark that for any concrete regularity structure $\mathscr{T} = \big((T^+,\Delta^+), (T,\Delta)\big)$, then 
$$
\mathscr{T}^+ := \big((T^+,\Delta^+), (T^+,\Delta^+)\big)
$$ 
is also a concrete regularity structure. For $g\in G^+$, set
\begin{equation}
\label{EqDefnHatPlus}
\widehat{g}^+\tau := (\textrm{Id}\otimes g)\Delta^+;
\end{equation}
this map sends $T^+$ into itself.

\medskip

\noindent \textsf{\textbf{Remark.}} \textit{For $g\in G^+$, the map $\widehat{g}$ is denoted by $\Gamma_g$ in Hairer's work \cite{Hai14}; we prefer the former Fourier-like notation}.

\medskip

We now come to the definition of the reference objects ${\sf\Pi}_x^{\sf g}\tau$ and ${\sf g}_{yx}(\sigma)$ used to give local descriptions of distributions and functions in a regularity structure setting, as in the introduction to this section. They come under the form of a model.

\bigskip

\hfill \textcolor{gray}{\textsf{\textbf{Models}}} 

\medskip

Recall $\beta_0 = \min A \in\bbR$. Given a function $\varphi$ on $\bbR^d$, and $x\in\bbR^d, 0<\lambda\leq 1$, set 
$$
\varphi_x^\lambda(\cdot) := \lambda^{-d}\varphi\big(\lambda^{-1}(\cdot-x)\big).   \vspace{0.1cm}
$$

\smallskip

\begin{defn*}
A \textbf{\textsf{model over a regularity structure}} $\mathscr{T}$ is a pair $({\sf g}, {\sf \Pi})$ of maps
$$
{\sf g} : \bbR^d \rightarrow G^+, \qquad {\sf \Pi} : T\rightarrow \mcC^{\beta_0}(\bbR^d)
$$
with the following properties.   \vspace{0.1cm}
\begin{itemize}
\item Set
$$
{\sf g}_{yx}:={\sf g}_y*{\sf g}_x^{-1}
$$
for each $x,y\in\bbR^d$. One has
\begin{align}\label{EqEstimateGammayx}
\|{\sf g}\|:=
\sup_{\tau\in\mcB^+}\sup_{x\in\bbR^d}|{\sf g}_x(\tau)|
+\sup_{\tau\in\mcB^+}\sup_{x,y\in\bbR^d}\frac{|{\sf g}_{yx}(\tau)|}{|y-x|^{|\tau|}}<\infty.
\end{align}
\item The map ${\sf\Pi}$ is linear. Set
$$
{\sf\Pi}_x^{\sf g}:=({\sf\Pi}\otimes{\sf g}_x^{-1})\Delta
$$
for each $x\in\bbR^d$. Fix $r > |\beta_0\wedge 0|$. One has
\begin{align}\label{EqestimatePix}
\|{\sf\Pi}\|^{{\sf g}}:=
\sup_{\sigma\in\mcB}\|{\sf\Pi}\sigma\|_{\mcC^{\beta_0}}
+\sup_{\sigma\in\mcB}\sup_{\varphi, 0<\lambda\le1, x\in\bbR^d}\frac{\big|\big\langle{\sf\Pi}_x^{\sf g}\sigma , \varphi_x^\lambda \big\rangle\big|}{\lambda^{|\sigma|}}<\infty,
\end{align}
where $\varphi$ runs over all functions $\varphi\in \mcC^r(\bbR^d)$, with associated norm no greater than $1$ and support in the unit ball.
\end{itemize}
\end{defn*}

\medskip

In Hairer's original work \cite{Hai14}, the notations $\Pi_x$ and $\Gamma_{yx}$ are used instead of ${\sf\Pi}_x^{\sf g}$ and $\widehat{{\sf g}_{yx}}$, respectively. In \eqref{EqEstimateGammayx} and \eqref{EqestimatePix}, we assume global bounds over $\bbR^d$, while Hairer only assumes in \cite{Hai14} the previous bounds in any compact subset of $\bbR^d$. In this paper, we work on the globally bounded case for simplicity. Our result may be extended into the locally bounded case using the weighted norms $\|f\|_{L_w^\infty} = \sup_{x\in\bbR^d} w^{-1}(x)|f(x)|$ instead of $\|f\|_{L^\infty}$.

For comparison, and given $a<0$, note that a distribution $\Theta$ on $\bbR^d$ is an element of $\mcC^a(\bbR^d)$ iff one has a bound 
$$
\big|\big\langle\Theta , \varphi^\lambda_x\big\rangle \big| \lesssim \lambda^a,
$$
for any $0<\lambda\le1$, uniformly in $x\in\bbR^d$ and $\varphi\in\mcC^r(\bbR^d)$, of unit norm in that space and support in the unit ball, for $r=|\lfloor a\rfloor|$. We stress that ${\sf\Pi}\tau$ is only an element of $\mcC^{\beta_0}(\bbR^d)$; identity \eqref{EqestimatePix} conveys the idea that ${\sf\Pi}^{\sf g}_x\tau$ behaves \textit{at point $x$} like an element of $\mcC^{|\tau|}(\bbR^d)$. Emphasize that ${\sf g}$ acts on $T^+$, while $\sf \Pi$ acts on $T$, and \textit{note that ${\sf g}$ plays on $T^+$ the same role as ${\sf \Pi}$ on $T$}; For $\tau\in T^+$ and $\sigma\in T$, one has
\begin{equation}
\label{EqParallel}
{\sf g}_{yx}(\tau) = \big({\sf g}(\cdot)(y)\otimes {\sf g}_x^{-1}\big)\Delta^+\tau, \quad
({\sf \Pi}_x^{\sf g}\sigma)(y) = \big({\sf \Pi}(\cdot)(y)\otimes {\sf g}_x^{-1}\big)\Delta\sigma,
\end{equation}
in a distributional sense for the latter. Note also the fundamental relation
\begin{equation}   \textcolor{red}{\label{EqTransitionRelation}}
{\sf\Pi}_y^{\sf g} = {\sf\Pi}_x^{\sf g}\circ\widehat{{\sf g}_{xy}},
\end{equation}
for all $x, y\in\bbR^d$; it comes from the comodule property \eqref{EqRightComodule}. The map ${\sf \Pi}$ can be recovered from each map ${\sf\Pi}_x^{\sf g}$, as we have 
\begin{equation}
\label{EqPiPix}
{\sf \Pi} = ({\sf\Pi}_x^{\sf g}\otimes {\sf g}_x)\Delta,
\end{equation}
as a consequence of the comodule property \eqref{EqRightComodule} 
\begin{equation*}
\begin{split}
({\sf\Pi}_x^{\sf g}\otimes {\sf g}_x)\Delta &= \big({\sf \Pi}\otimes {\sf g}_x^{-1}\otimes {\sf g}_x\big)(\Delta\otimes \textrm{Id})\Delta   \\
&= \big({\sf \Pi}\otimes {\sf g}_x^{-1}\otimes {\sf g}_x\big)(\textrm{Id}\otimes \Delta^+)\Delta   \\
&= ({\sf \Pi} \otimes {\bf 1}')\Delta = {\sf \Pi}.
\end{split}
\end{equation*}

\medskip

\noindent \textbf{\textsf{Examples.}} \textbf{\textsf{1. Bounded polynomials structure.}}  For any smooth function $f$ on $\bbR^d$, and $r>0$, the Taylor expansion property
$$
f(y)-\sum_{|k|<r}\frac{\partial^kf(x)}{k!}(y-x)^k=O(|y-x|^r).
$$
is usually lifted to a modelled distribution
$$
\bsf(x) := \sum_{|k|<r}\frac{\partial^kf(x)}{k!}\bsX^k,
$$
over the canonical polynomial regularity structure, under the model $({\sf \Pi}\bsX^k)(x)=x^k$ and ${\sf g}_x(\bsX^k)=x^k$. Since they are not bounded functions, we modify this expansion by using smooth and bounded functions behaving like polynomials in local sets. The following elementary claim is proved in Appendix \ref{SectionBoundedPolynomials}. 

\medskip

\begin{prop}\label{app: prop PTE}
There exists a finite set $E$, an open covering $\{U_e\}_{e\in E}$ of $\bbR^d$, and a family $\big\{\phi_e,\{x\mapsto x_e^i\}_{i=1}^d\big\}_{e\in E}$ of functions enjoying the following properties.
\begin{enumerate}
\item[\sf (a)] The functions $\phi_e : \bbR^d\to[0,\infty)$, belong to $C^\infty_b(\bbR^d)$, $\phi_e(x)=0$ for any $x\in U_e^c$, and $\sum_{e\in E}\phi_e(x)=1$ for any $x\in\bbR^d$.
\item[\sf (b)] The functions $x\mapsto x_e^i$, belong to $C_b^\infty(\bbR^d)$, and $y_e^i-x_e^i=y^i-x^i$ for any $x,y$ on the connected component of $U_e$.
\item[\sf (c)] For any $f\in C^\infty_b(\bbR^d)$ and $r>0$, we have
\begin{align}\label{app: patched Taylor expansion}
\Big|f(y)-\sum_{e \in E}\sum_{|k|<r}\frac{\partial^k(\phi_e f)(x)}{k!}(y_e - x_e)^k\Big|\lesssim B_r(f) \, |y-x|^r,
\end{align}
where $B_r(f):=\|f\|_{C_b^r}$, if $r\in\bbN$, or $B_r(f):=\|f\|_{\mcC^r}$, if $r\in(0,\infty)\setminus\bbN$.
\end{enumerate}
\end{prop}

\medskip

We lift expansion \eqref{app: patched Taylor expansion} to an appropriate regularity structure as follows. Let 
$$
X := \{X_e^i\}_{e\in E,1\leq i\leq d}
$$
be a family of symbols, and let $T^+(X)$ be the commutative free algebra with unit $\unit$, generated by these symbols. We define a coproduct $\Delta^+:T^+(X)\to T^+(X)\otimes T^+(X)$ by
$$
\Delta^+{\bf1}={\bf1}\otimes{\bf1},\quad
\Delta^+X_e^i = X_e^i\otimes{\bf1} + {\bf1}\otimes X_e^i,
$$
which turns $T^+(X)$ into a Hopf algebra. By defining the homogeneity $|\cdot|$ by $|X_e^i|=1$, we have the graded Hopf algebra $T^+(X)$.
Let $T(X)$ be the subspace spanned by the bounded polynomials $\{\bsX_e^k\}_{e\in E,k\in\bbN^d}$, where\vspace{-0.2cm}
$$
\bsX_e^k := \prod_{i=1}^d \big(X_e^i\big)^{k_i},\quad k=(k_i)_{i=1}^d\in\bbN^d.
$$
Denote by 
$$
\Delta:T(X)\to T(X)\otimes T^+(X)
$$
the restriction of $\Delta^+$ to $T(X)$, which turns $T(X)$ into a right comodule over $T^+(X)$. By definition, we have the concrete regularity structure $\mathscr{T}(X) := \big(T^+(X), T(X)\big)$. The canonical model $\sf (g, \Pi)$ on $\mathscr{T}(X)$ is defined by
\begin{equation}
\label{EqModelsBoundedPolynomials}
{\sf g}_x(\bsX_e^k) = ({\sf \Pi}\bsX_e^k)(x) = x_e^k.
\end{equation}
The following elementary result, proved in Appendix \ref{SectionBoundedPolynomials} provides the canonical lift of a smooth function to this bounded polynomials regularity structure. See the paragraph on modelled distributions for the definition of $\mcD^r(\mathscr{T}(X), {\sf g})$ and the associated norm $\trino{\cdot}_{\mcD^r}$.

\medskip

\begin{prop}
\label{PropPolynomialLift}
For any given $f\in C_b^\infty(\bbR^d)$ and $r>0$, define the $T(X)$-valued function
$$
\bsf(x) := \sum_{e\in E}\sum_{|k|<r} \frac{\partial^k(\phi_e f)(x)}{k!} \bsX_e^k,  \quad x\in\Euc^d.
$$
Then $\bsf\in\mcD^r(\mathscr{T}(X),{\sf g})$, and $\trino{\bsf}_{\mcD^r}\lesssim B_r(f)$.
\end{prop}

\medskip

\textbf{\textsf{2. Canonical model on $\mathscr{T}^+$.}} As another example of model over some regularity structure, consider the regularity structure $\mathscr{T}^+$ associated with any regularity structure $\mathscr{T}$, and assume we are given a function ${\sf g} : \bbR^d\mapsto G^+$ that satisfies estimate \eqref{EqEstimateGammayx}. For $\tau\in T^+$, set 
\begin{equation}
\label{EqCanonicalModel}
{\sf\Pi}^{\sf g}\tau(x):={\sf g}_x(\tau).
\end{equation}
Estimate \eqref{EqestimatePix} holds as a consequence of \eqref{EqEstimateGammayx}, so $({\sf g},{\sf\Pi}^{\sf g})$ is a model on $\mathscr{T}^+=(T^+,T^+)$. This justifies to say simply that ${\sf g}$ is a model on $\mathscr{T}^+$.   \hfill $\ovoid$

\bigskip

Equation \eqref{EqPiPix} giving $\sf \Pi$ in terms of ${\sf\Pi}^{\sf g}_x$ and ${\sf g}_x$ writes explicitly 
$$
{\sf \Pi}\tau = \sum_{\sigma\leq \tau} {\sf\Pi}^{\sf g}_x(\sigma) {\sf g}_x(\tau/\sigma),\vspace{-0.15cm}
$$
for $\tau\in\mcB$, that is\vspace{-0.1cm}
\begin{equation}
\label{EqPixTau}
{\sf\Pi}^{\sf g}_x\tau = {\sf \Pi}\tau - \sum_{\sigma<\tau} {\sf g}_x(\tau/\sigma){\sf\Pi}^{\sf g}_x\sigma.
\end{equation}
Furthermore expanding ${\sf\Pi}_x^{\sf g}\sigma$, one has
$$
{\sf\Pi}^{\sf g}_x\tau = {\sf \Pi}\tau - \sum_{\sigma_1<\tau} {\sf g}_x(\tau/\sigma_1){\sf \Pi}\sigma_1
+\sum_{\sigma_2<\sigma_1<\tau}{\sf g}_x(\tau/\sigma_1){\sf g}_x(\sigma_1/\sigma_2){\sf\Pi}^{\sf g}_x\sigma_2.
$$
Iterating this expansion gives a representation of ${\sf\Pi}^{\sf g}_x$ in terms of ${\sf g}_x$ and ${\sf \Pi}$
\begin{equation}
\label{EqPixIteratePi}
{\sf\Pi}^{\sf g}_x\tau = {\sf \Pi}\tau - \sum_{n\geq 1}(-1)^{n-1}\sum_{\sigma_n<\cdots<\sigma_1<\tau} {\sf g}_x(\tau/\sigma_1)\cdots {\sf g}_x(\sigma_{n-1}/\sigma_n)\, {\sf \Pi}\sigma_n;
\end{equation}
the sum is finite. Similarly, since ${\sf g}_y = {\sf g}_{yx}*{\sf g}_x$, by definition, Lemma \ref{LemmaCoprodQuotient} provides for any $\sigma\leq^{(+)}\tau\in\mcB^{(+)}$ the relation
$$
{\sf g}_{yx}\big(\tau/^{(+)}\sigma\big) = {\sf g}_y\big(\tau/^{(+)}\sigma\big) - {\sf g}_x\big(\tau/^{(+)}\sigma\big) - \sum_{\sigma<^{(+)}\sigma_1<^{(+)}\tau} {\sf g}_x\big(\tau/^{(+)}\sigma_1\big) \, {\sf g}_{yx}\big(\sigma_1/^{(+)}\sigma\big).
$$
A repeated expansion then gives a representation of ${\sf g}_{yx}(\tau/^{(+)}\sigma)$ in terms of ${\sf g}_y$ and ${\sf g}_x$
\begin{equation}
\label{EqGammayxIterateG}
\begin{split}
&{\sf g}_{yx}\big(\tau/^{(+)}\sigma\big) = {\sf g}_y\big(\tau/^{(+)}\sigma\big) - {\sf g}_x\big(\tau/^{(+)}\sigma\big)   \\
&- \sum_{n\geq 1} (-1)^{n-1} \hspace{-0.35cm}\sum_{\sigma<^{(+)}\sigma_n<^{(+)}\cdots<^{(+)}\tau} {\sf g}_x\big(\tau/^{(+)}\sigma_1\big)\cdots {\sf g}_x\big(\sigma_{n-1}/^{(+)}\sigma_n\big)\,\Big({\sf g}_y\big(\sigma_n/^{(+)}\sigma\big) - {\sf g}_x\big(\sigma_n/^{(+)}\sigma\big)\Big).
\end{split} \end{equation}

\bigskip

\hfill \textcolor{gray}{\textsf{\textbf{Modelled distributions}}}

\medskip

Recall notation \eqref{EqNotationDirectSum} for the notation $\|h\|_\alpha$ for $\alpha\in A$ and $h\in T$.

\medskip

\begin{defn*}
Let ${\sf g} : \bbR^d\rightarrow G^+$ satisfy \eqref{EqEstimateGammayx}. Fix a regularity exponent $\gamma\in\bbR$. One defines the space $\mcD^\gamma(\mathscr{T}, {\sf g})$ of \textbf{\textsf{distributions modelled on the regularity structure $\mathscr{T}$, with transition ${\sf g}$}}, as the space of functions $\bsf : \bbR^d\rightarrow T_{<\gamma}$ such that
\begin{equation*}
\begin{split}
&\brarb{\bsf}_{\mcD^\gamma} := \max_{\beta<\gamma}\, \sup_{x\in\bbR^d}\, \big\| \bsf(x) \big\|_{\beta} < \infty,   \\
&\|\bsf\|_{\mcD^\gamma} := \max_{\beta<\gamma}\sup_{x,y\in\bbR^d}\frac{\big\| \bsf(y) - \widehat{{\sf g}_{yx}}\bsf(x) \big\|_{\beta}}{|y-x|^{\gamma-\beta}}<\infty.
\end{split}
\end{equation*}
Set $\trino{\bsf}_{\mcD^\gamma} := \brarb{\bsf}_{\mcD^\gamma} + \|\bsf\|_{\mcD^\gamma}$.
\end{defn*}

\medskip

For a basis element $\sigma\in\mcB$, and an arbitrary element $h$ in $T$, denote by $h^\sigma$ its component along the $\sigma$ direction. For a modelled distribution $\bsf(\cdot) = \sum_{\sigma\in\mcB} f^\sigma(\cdot) \, \sigma$ in $\mcD^\gamma(\mathscr{T}, {\sf g})$, and $\sigma_0\in\mcB$, we have
\begin{equation}
\label{EqIncrementMdF}
\Big(\bsf(y) - \widehat{{\sf g}_{yx}}\bsf(x)\Big)^{\sigma_0} = f^{\sigma_0}(y) - f^{\sigma_0}(x) - \sum_{\tau > \sigma_0} {\sf g}_{yx}(\tau/\sigma_0) \, f^\tau(x).
\end{equation}
As an example, given a basis element $\tau\in \mcB$, set
\begin{equation}
\label{EqDefnMdTau}
\bsh_\tau(x) := \sum_{\sigma<\tau}{\sf g}_x(\tau/\sigma)\sigma.
\end{equation}
Then, it follows from identity \eqref{EqCoprodQuotient} giving $\Delta^+(\tau/\sigma)$, in Lemma \ref{LemmaCoprodQuotient}, that
\begin{equation*}
\begin{split}
\widehat{{\sf g}_{yx}}\bsh_\tau(x) &= \sum_{\eta\leq \sigma<\tau} {\sf g}_{yx}(\sigma/\eta){\sf g}_x(\tau/\sigma)\eta = \sum_{\eta<\tau} \big({\sf g}_y(\tau/\eta) - {\sf g}_{yx}(\tau/\eta)\big)\eta   \\
&= \bsh_\tau(y) - \sum_{\eta<\tau}{\sf g}_{yx}(\tau/\eta)\eta.
\end{split}
\end{equation*}
The size estimate $\big|{\sf g}_{yx}(\tau/\eta)\big| \lesssim |y-x|^{|\tau|-|\eta|}$, then shows that $\bsh_\tau$ is a modelled distribution in $\mcD^{|\tau|}(\mathscr{T}, {\sf g})$. Here is another example.

\medskip

\begin{lem}   \label{LemmaExModelledFunction}
Let $\bsf = \sum_{\sigma\in\mcB}f^\sigma(\cdot)\,\sigma$, be an element of $\mcD^\gamma(\mathscr{T},{\sf g})$. Then, for each $\tau\in\mcB$, the $T^+$-valued function 
$$
\bsf/\tau := \sum_{\sigma\geq\tau}f^\sigma(\cdot)\,\sigma/\tau.
$$
is an element of $\mcD^{\gamma-\vert\tau\vert}(\mathscr{T}^+,{\sf g})$.
\end{lem}

\medskip

\begin{Dem}
This comes from the identity 
$$
(\bsf/\tau)(y) - \widehat{{\sf g}_{yx}}^+(\bsf/\tau)(x) = \sum_{\sigma\geq\tau}\Big(f^\sigma(y) - \sum_{\mu\geq\sigma} f^\mu(x)\,{\sf g}_{yx}(\mu/\sigma)\Big)\,\sigma/\tau,
$$
and the fact that $\bsf$ is a modelled distribution.
\end{Dem}

\medskip

Recall $\beta_0 = \min A$, and fix $r>|\beta_0\wedge0|$.

\medskip

\begin{thm}   \label{ThmReconstructionRS}
\textbf{\textsf{(Hairer's reconstruction theorem)}}
Let $({\sf g},{\sf \Pi})$ be a model over $\mathscr{T}$. Fix a regularity exponent $\gamma\in\bbR$. There exists a linear continuous operator 
$$
\textsf{\textbf{R}} : \mcD^\gamma(\mathscr{T}, {\sf g}) \rightarrow \mcC^{\beta_0}(\bbR^d)
$$
satisfying the property
\begin{equation}
\label{EqReconstructionCondition}
\Big|\big\langle \textsf{\textbf{R}}\bsf - {\sf\Pi}_x^{\sf g}\bsf(x) , \varphi^\lambda_x\big\rangle\Big|  \lesssim  \|{\sf \Pi}\|^{\sf g}\, \big\|\bsf\big\|_{\mcD^\gamma} \, \lambda^\gamma,
\end{equation}
uniformly in $\bsf\in\mcD^\gamma(\mathscr{T}, {\sf g}), \varphi\in\mcC^r(\bbR^d)$ with unit norm and support in the unit ball, $x\in\bbR^d$ and $0<\lambda\leq 1$. Such an operator is unique if the exponent $\gamma$ is positive.
\end{thm}

\medskip

A distribution satisfying identity \eqref{EqReconstructionCondition} is called a \textit{reconstruction of the modelled distribution} $\bsf$. See Theorem 3.10 in Hairer' seminal work \cite{Hai14}. We provide in Theorem \ref{ThmReconstructionPC} below an explicit representation for the reconstruction operator $\textsf{\textbf{R}}$ building on paracontrolled calculus. Notice from the definition of ${\sf\Pi}^{\sf g}_x$ that the constraint $\big|\big\langle{\sf\Pi}^{\sf g}_x\tau , \varphi_x^\lambda \big\rangle\big|  \lesssim \lambda^{|\tau|}$, that needs to be satisfied by a model, is equivalent to the estimate
$$
\Big|\big\langle {\sf \Pi}\tau - \sum_{\sigma<\tau} {\sf g}_x(\tau/\sigma){\sf\Pi}^{\sf g}_x\sigma , \varphi^\lambda_x\big\rangle\Big|  \lesssim  \lambda^{|\tau|}.
$$
This means that ${\sf \Pi}\tau$, with $\tau\in \mcB$, is a reconstruction of the modelled distribution $\bsh_\tau\in\mcD^{|\tau|}(\mathscr{T},{\sf g})$ defined above in \eqref{EqDefnMdTau}. Recall that uniqueness in the reconstruction theorem implies that if $\bsf$ takes values in a function-like sector of $T$, then $\textsf{\textbf{R}}\bsf = {\bf 1}'(\bsf)$ -- see e.g. Proposition 3.28 in Section 3.4 of \cite{Hai14}.

\bigskip

\section{Explicit formula for the reconstruction operator}
\label{SubsectionExplicitReconstruction}

We prove Theorem \ref{ThmMain1} giving an explicit description of the reconstruction operator in this section.


\subsection{From Taylor local description to global paracontrolled representation}

We describe here some simple properties of a natural two-parameter extension of the elementary paraproduct built from Littlewood-Paley blocks, and refer the reader to Appendix \ref{AppendixParapdcts} for background on Littlewood-Paley decomposition. The notations $\Delta_i$ and $Q_i$ for the $i^\textrm{th}$ Littlewood-Paley block and its kernel are recalled in Appendix \ref{AppendixParapdcts}. For $j\geq 1$, define the operator $S_j := \sum_{i\le j-2} \Delta_i$, and its smooth kernel $P_j := \sum_{i\le j-2} Q_i$. The H\"older spaces $\mcC^\alpha(\bbR^d)$ are defined as Besov spaces $B^\alpha_{\infty\infty}(\bbR^d)$ -- see Appendix \ref{AppendixParapdcts}.

\medskip

For a two-variable real-valued distribution $\Lambda$ on $\bbR^d\times\bbR^d$, and $j\geq 1$, set
$$
(\bfQ_j\Lambda)(x) := \iint P_j(x-y) Q_j(x-z)\,\Lambda(y,z)dydz;
$$
we abuse notations using the integral notation. Set 
$$
\bfP\Lambda  :=  \sum_{j\geq 1} \bfQ_j\Lambda.
$$
We often write 
$$
\bfP \Lambda = \bfP_{y,z}\big(\Lambda(y,z)\big)
$$ 
in order to display the integrated variables. With that notation, we have the consistency relation
$$
{\sf P}_fg = \bfP_{y,z}\big(f(y)g(z)\big), \quad f,g\in L^\infty,
$$
between the paraproduct operator ${\sf P}$ and its two-parameter extension. For $\alpha>0$, and a measurable real-valued function $F$ on $\bbR^d\times\bbR^d$, set 
$$
\trino{F}_{\mcC_2^\alpha} := \sup_{y,z\in\bbR^d}\frac{|F(y,z)|}{|y-z|^\alpha}.
$$

\medskip

\begin{prop}
\label{PropParaprodEstimates}
\begin{enumerate}
   \item[\textsf{\emph{(a)}}] Let $\Lambda$ be a real-valued distribution on $\bbR^d\times\bbR^d$. If $\|\bfQ_j\Lambda\|_{L^\infty}\lesssim2^{-j\alpha}$, for all $j\ge1$, for some $\alpha\in\bbR$, then ${\bfP}\Lambda\in\mcC^\alpha(\bbR^d)$, and
   $$
   \|{\bfP} \Lambda\|_{\mcC^\alpha}\lesssim\sup_{j\ge1}2^{j\alpha}\|\bfQ_j\Lambda\|_{L^\infty}.
   $$

   \item[\textsf{\emph{(b)}}] Let $\alpha>0$, and a real-valued measurable function $F$ on $\bbR^d\times\bbR^d$ be given, with $\trino{F}_{\mcC_2^\alpha}<\infty$. Then ${\bfP}F\in\mcC^\alpha(\bbR^d)$, and $\|{\bfP} F\|_{\mcC^\alpha}\lesssim\trino{F}_{\mcC_2^\alpha}$.
\end{enumerate}
\end{prop}

\medskip

\begin{Dem}
\textsf{(a)} Since $\mathscr{F} P_j$ is supported in $\big\{x\in\bbR^d;|x|<2^j\times\frac23\big\}$ and $\mathscr{F} Q_j$ is supported in $\big\{x\in\Euc^d;2^j\times\frac34<|x|<2^j\times\frac83\big\}$, the integral
\begin{align*}
\int Q_i(x-w)P_j(w-y)Q_j(w-z)dw
\end{align*}
vanishes if $|i-j|\ge5$. Hence $\Delta_i({\bfP} \Lambda)=\sum_{|i-j|\le4}\Delta_i(\bfQ_j\Lambda)$ and we have
$$
\|\Delta_i({\bfP} \Lambda)\|_{L^\infty}
\le\sum_{|i-j|\le4}\|\Delta_i(\bfQ_j\Lambda)\|_{L^\infty}
\lesssim\sum_{|i-j|\le4}\|\bfQ_j\Lambda\|_{L^\infty}
\lesssim\sum_{|i-j|\le4}2^{-\alpha j}
\lesssim2^{-\alpha i}.
$$

\textsf{(b)} It is sufficient to show that $\|\bfQ_jF\|_{L^\infty}\lesssim2^{-\alpha j}$ for all $j\ge2$. By the scaling properties $P_j(\cdot)=2^{(j-2)d}P_{2}(2^{j-2}\,\cdot)$ and $Q_j(\cdot)=2^{(j-2)d}Q_2(2^{j-2}\,\cdot)$, we have
\begin{align*}
|\bfQ_jF(x)|
&\lesssim \iint |P_j(x-y)Q_j(x-z)||y-z|^\alpha dydz\\
&= 2^{-\alpha (j-2)}\iint |P_2(2^{j-2}x-y)Q_2(2^{j-2}x-z)||y-z|^\alpha dydz\\
&= 2^{-\alpha (j-2)}\iint |P_2(y)Q_2(z)||y-z|^\alpha dydz\lesssim 2^{-\alpha j}.
\end{align*}
\end{Dem}

\medskip

The next proposition is the key step to the representation of the reconstruction operator given in Theorem 6.10 of \cite{GIP15}. We state it and prove it here under a slightly more general form. See the proofs of Lemma 6.8, Lemma 6.9 and Theorem 6.10 therein.

\medskip

\begin{prop}
\label{PropReconstructionPC}
Let $\gamma\in\bbR$ and $\beta_0\in\bbR$ be given together with a family $\Lambda_x$ of distributions on $\bbR^d$, indexed by $x\in\bbR^d$. Assume that one has
$$
\sup_{x\in\bbR^d}\|\Lambda_x\|_{\mcC^{\beta_0}}< \infty
$$
and one can decompose $(\Lambda_y - \Lambda_x)$ under the form
   \begin{equation}
   \label{EqStructureWorkHorse}
   \Lambda_y - \Lambda_x = \sum_{\ell=1}^N c^\ell_{yx}\,\Theta_x^\ell   
   \end{equation}
for $N$ finite, $\bbR^d$-indexed distributions $\Theta_x^\ell$, and real-valued coefficients $c_{yx}^\ell$ depending measurably on $x$ and $y$, such that 
$$
\sup_{x\in\bbR^d} \, \sup_{j\ge-1}\, 2^{j\beta_\ell} \big| \big\langle \Theta^\ell_x, P_j(x-\cdot)\big\rangle \big| < \infty, \qquad \trino{c^\ell}_{\mcC_2^{\gamma-\beta_\ell}} < \infty,
$$
for regularity exponents $\beta_\ell < \gamma$, for all $1\leq \ell\leq N$. Denote ${\bfP}(\Lambda)={\bfP}_{y,z}(\Lambda_y(z))$ below.  \vspace{0.15cm}

\begin{itemize}
   \item If $\gamma>0$, then there exists a unique function $f_\Lambda\in\mcC^\gamma(\bbR^d)$ such that 
   \begin{equation}
   \label{EqReconstructionPC1}
   \Big|\Big\langle\big\{{\bfP}(\Lambda) - f_\Lambda\big\} - \Lambda_x, P_i(x-\cdot)\Big\rangle\Big| \lesssim 2^{-i\gamma},   
   \end{equation}
   uniformly in $x\in\bbR^d$.   \vspace{0.1cm}
   
   \item If $\gamma<0$, then we have
   \begin{equation}
   \label{EqReconstructionPC2}
   \big|\big\langle{\bfP}(\Lambda) - \Lambda_x, P_i(x-\cdot)\big\rangle\big| \lesssim 2^{-i\gamma},   
   \end{equation}
   uniformly in $x\in\bbR^d$.
\end{itemize}
\end{prop}

\medskip

\begin{Dem}
\noindent \textsf{\textbf{(i)}}
We prove that one has
\begin{equation}
\label{EqEstimatePCReconstruction}
\Big|\Delta_j\big({\bfP}(\Lambda)-\Lambda_x\big)(x)\Big| \lesssim 2^{-j\gamma},
\end{equation}
uniformly in $x\in\bbR^d$. We write for that purpose
\begin{equation*}
\begin{split}
{\bfP}(\Lambda)(y) - \Lambda_x(y) &= 
\sum_{j\geq -1} \iint P_j(y-u)Q_j(y-v)\big(\Lambda_u(v)-\Lambda_x(v)\big) \, dudv - \mathscr{S}(\Lambda_x)   \\
&= \sum_{1\leq \ell\leq N}\sum_{j\geq -1} \iint P_j(y-u)Q_j(y-v)c_{ux}^\ell \Theta^\ell_x(v) \, dudv - \mathscr{S}(\Lambda_x).
\end{split}
\end{equation*}
Here the operator $\mathscr{S}$ is defined by
\begin{equation}
\label{EqDefnMscrS}
\mathscr{S}f := f- {\sf P}_1f = f-{\bfP}_{y,z}\big(1(y)f(z)\big)
\end{equation}
for any $f\in\mcS'(\bbR^d)$. This is a smooth function that depends continuously on $f$; if $f\in\mcC^\alpha(\bbR^d)$ with $\alpha\in\bbR$, then for any $r>0$,
$$
\|\mathscr{S}f\|_{\mcC^r}\lesssim\|f\|_{\mcC^\alpha}.
$$
Hence we have for any $i\geq1$,
$$
\big|\Delta_i\big({\bfP}(\Lambda) - \Lambda_x\big)(x)\big| \leq \sum_{|j-i|\leq 4} \sum_\ell\iint \big|Q_i(x-y)P_j(y-u)\big|\,|u-x|^{\gamma-\beta_\ell}2^{-i\beta_\ell} \, dudy + o(2^{-i\gamma}).
$$
Then \eqref{EqEstimatePCReconstruction} follows from elementary estimates and the bounds
\begin{equation}
\label{EqEstimateLPProjector}
\int_{\bbR^d} |P_j|(x) \, |x|^r \, dx \lesssim 2^{-jr},  \qquad \int_{\bbR^d} |Q_j|(x) \, |x|^r \, dx \lesssim 2^{-jr},
\end{equation}
that holds for any positive exponent $r$.

\medskip

\noindent \textsf{\textbf{(ii)}} If $\gamma>0$, estimate \eqref{EqEstimatePCReconstruction} implies that the sum
$$
f_\Lambda(x) := \sum_{j\geq -1}\Delta_j\big({\bfP}(\Lambda)-\Lambda_x\big)(x),
$$
defines an element $f_\Lambda$ of $\mcC^\gamma(\bbR^d)$; this is proved in point \textsf{\textbf{(iii)}} below. Then we have, for any $x\in\bbR^d$,
\begin{align*}
&\big|\big\langle{\bfP}(\Lambda)-f_\Lambda - \Lambda_x, P_i(x-\cdot)\big\rangle\big| = \left| \sum_{j\le i-2}\Delta_j\big({\bfP}(\Lambda)-\Lambda_x\big) - S_i(f_\Lambda)(x)\right|   \\
&=\left| f_\Lambda(x) - \sum_{j> i-2}\Delta_j\big({\bfP}(\Lambda)-\Lambda_x\big) - S_i(f_\Lambda)(x)\right|  \\
&\lesssim \sum_{j>i-2}\big|\Delta_j(f_\Lambda)(x)\big| + \sum_{j>i-2}\big|\Delta_j\big({\bfP}(\Lambda)-\Lambda_x\big)(x)\big|\lesssim 2^{-i\gamma}.
\end{align*}
Uniqueness of $f_\Lambda$ follows from the fact that $P_i$ converges to a Dirac mass at $0$, so if $f_\Lambda'$ were another $\mcC^\gamma$ function satisfying estimate \eqref{EqReconstructionPC1}, one would have
$$
\big|\big\langle f_\Lambda-f'_\Lambda, P_i(x-\cdot)\big\rangle\big| \lesssim 2^{-i\gamma},
$$
uniformly in $x$, for all $i\geq -1$, giving indeed $f_\Lambda'=f_\Lambda$.

\ssk

If $\gamma<0$, we directly have from \eqref{EqEstimatePCReconstruction} that 
$$
\big|\big\langle{\bfP}(\Lambda) - \Lambda_x, P_i(x-\cdot)\big\rangle\big| \lesssim \sum_{j\le i-2}\big|\big\langle{\bfP}(\Lambda) - \Lambda_x, Q_j(x,\cdot)\big\rangle\big| \lesssim \sum_{j\le i-2}2^{-j\gamma}\lesssim2^{-i\gamma}.
$$

\medskip

\noindent \textsf{\textbf{(iii)}} We follow the argument in Section 6 of \cite{GIP15}. We decompose $f_\Lambda=f^{\le j+1}_\Lambda+f^{>j+1}_\Lambda$, where
$$
f_\Lambda^{\le j+1}(x) := \sum_{i\le j+1}\Delta_i\big({\bfP}(\Lambda)-\Lambda_x\big)(x).
$$
We consider $\Delta_jf_\Lambda=\Delta_jf_\Lambda^{\le j+1}+\Delta_jf_\Lambda^{>j+1}$. For the second term, by the estimate \eqref{EqEstimatePCReconstruction} one has
$$
\big\|\Delta_jf_\Lambda^{>j+1}\big\|_{L^\infty}\lesssim\sum_{i>j+1}2^{-i\gamma}\lesssim2^{-j\gamma}.
$$
For the first term, since $Q_j*Q_{\le j+1}=Q_j$,  one has
\begin{align*}
\Delta_jf_\Lambda^{\le j+1}(y)
&=\int Q_j(y-x)Q_{\le j+1}\big({\bfP}(\Lambda)-\Lambda_x\big)(x)dx   \\
&=\int Q_j(y-x)Q_{\le j+1}\Big({\bfP}(\Lambda)-\Lambda_y+\sum_{\ell} c^\ell_{yx}\Theta_x^\ell\Big)(x)dx   \\
&=Q_j\big({\bfP}(\Lambda)-\Lambda_y\big)(y)+\sum_{\ell}\int Q_j(y-x)\,c_{yx}^\ell\, \big(Q_{\le j+1}\Theta_x^\ell\big)(x)\,dx.
\end{align*}
The first term is estimated by \eqref{EqEstimatePCReconstruction}. The second term is bounded by $2^{-j\gamma}$ by assumption. In the end, we have
$$
\big\|\Delta_jf_\Lambda^{\le j+1}\big\|_{L^\infty} \lesssim 2^{-j\gamma}.
$$
\end{Dem}

\smallskip

If $\Lambda_x$ stand for ${\sf\Pi}_x^{\sf g}\bsf(x)$, for a modelled distribution $\bsf\in\mcD^\gamma(\mathscr{T},{\sf g})$ and a model $\sf (g,\Pi)$, one has 
$$
\Lambda_y - \Lambda_x = \sum_{\sigma\in\mcB}\Big(\bsf(y) - \widehat{{\sf g}_{xy}}\bsf(y)\Big)^\sigma\,{\sf\Pi}^{\sf g}_x\sigma,
$$
and $\Lambda$ satisfies the assumptions of Proposition \ref{PropReconstructionPC}, from condition \eqref{EqestimatePix} on models  and the definition of a modelled distribution. As in Lemma 6.3 of \cite{GIP15}, we can extend the condition \eqref{EqestimatePix} for any rapidly decreasing smooth functions $\varphi$. Identities \eqref{EqReconstructionPC1} and \eqref{EqReconstructionPC2} are equivalent to saying that ${\bfP}(\Lambda)-f_\Lambda{\bf 1}_{\gamma>0}$ is a reconstruction of $\bsf$ -- see Lemma 6.6 of \cite{GIP15}. This is the content of Theorem 6.10 in \cite{GIP15}.

\medskip

We prove in Theorem \ref{ThmReconstructionPC} below that ${\bfP}_{y,z}\big(({\sf\Pi}_y^{\sf g}\bsf(y))(z)\big)$ has an explicit form, up to some remainder in $\mcC^\gamma(\bbR^d)$. The mechanism at work in the proof of this fact lies in Proposition \ref{PropWorkingHorse}. Following \cite{BB16}, set 
$$
{\sf R}^\circ(a,b,c) := {\sf P}_a({\sf P}_bc)-{\sf P}_{ab}c.
$$
It was proved in Appendix C1 of \cite{BB16} that the map ${\sf R}^\circ$ is continuous from $L^\infty(\bbR^d)\times\mcC^{r_1}(\bbR^d)\times\mcC^{r_2}(\bbR^d)$ into $\mcC^{r_1+r_2}(\bbR^d)$, for $r_1\in(0,1)$ and $r_2$ any regularity exponent in $\bbR$. The next proposition provides a refined continuity result for the operator $\sf R$.

\medskip

\begin{prop}
\label{PropWorkingHorse}
Pick a positive regularity exponent $\alpha$. Assume we are given a function $f\in L^\infty(\bbR^d)$ and a finite family $(a_n,b_n)_{1\leq n\leq N}$ of elements of $L^\infty(\bbR^d)\times L^\infty(\bbR^d)$ such that one has
\begin{equation}
\label{EqWorkingHorse}
f(y)-f(x) = \sum_{n=1}^N a_n(x)\big(b_n(y)-b_n(x)\big) + f^\sharp_{yx},
\end{equation}
for any $x,y\in\bbR^d$, for a two-parameter remainder $f^\sharp$ with finite $\alpha$-H\"older norm $\trino{f^\sharp}_{\mcC_2^\alpha}<\infty$. Then, for any regularity exponent $\beta\in\bbR$ and $g\in\mcC^\beta(\bbR^d)$, we have
$$
\sum_{n=1}^N{\sf R}^\circ(a_n,b_n,g) \in \mcC^{\alpha+\beta}(\bbR^d).
$$
\end{prop}

\medskip

\begin{Dem}
Recall from equation \eqref{EqDefnMscrS} the definition of the smooth function $\mathscr{S}g$, for any $g\in\mcC^\beta(\bbR^d)$, with $\beta\in\bbR$, and note the identity
$$
{\sf R}^\circ(a,b,c) = {\sf P}\Big(a(x)\big({\sf P}_{b-b(x)}c\big)(y)\Big) - {\sf P}_{ab}(\mathscr{S}c).
$$
Applying the two-parameter $\bfP$-operator to identity \eqref{EqWorkingHorse}, we see that 
\begin{equation*} \begin{split}
\sum_{n=1}^N {\sf R}^\circ(a_n,b_n,g) &= {\sf R}^\circ(1,f,g) + {\sf P}_f(\mathscr{S}g) - \sum_{n=1}^N{\sf P}_{a_nb_n}(\mathscr{S}g) - {\bfP}_{x,y}\Big(({\sf P}_{f^\sharp_{\cdot x}}g)(y)\Big)   \\
&= -\mathscr{S}({\sf P}_fg) + {\sf P}_f(\mathscr{S}g) - \sum_{n=1}^N{\sf P}_{a_nb_n}(\mathscr{S}g) - {\bfP}_{x,y}\Big(({\sf P}_{f^\sharp_{\cdot x}}g)(y)\Big).
\end{split} \end{equation*}
The first three terms on the right hand side are smooth. To prove that the last term on the right hand side is an element of $\mcC^{\alpha+\beta}(\bbR^d)$, it is sufficient, from Proposition \ref{PropParaprodEstimates}, to see that 
$$
\Big|{\bf Q}_j\Big(\big({\sf P}_{f^\sharp_{\cdot x}}g\big)(y)\Big)\Big| \lesssim 2^{-j(\alpha+\beta)},
$$
for all $j\geq 1$. Recall for that purpose the bound \eqref{EqEstimateLPProjector}. Then we have for $\Big|{\bf Q}_j\Big(\big({\sf P}_{f^\sharp_{\cdot x}}g\big)(y)\Big)\Big|$ the upper bound
\begin{equation*} \begin{split}
&\sum_{i;|i-j|\leq 4}\left|\int P_j(z-x)Q_j(z-y) 
\left(\int P_i(y-u)Q_i(y-v) f^\sharp_{ux} g(v) dudv\right)
\,dxdy   \right|\\
&\lesssim \sum_{i;|i-j|\leq 4} \int \big|P_j(z-x)Q_j(z-y)\big| \,\left(\int \big|P_i(y-u)\big|\,|u-x|^\alpha\,du\right)\,|\Delta_ig(y)|\,dxdy   \\
&\lesssim \sum_{i;|i-j|\leq 4} 2^{-i\beta}\int \big|P_j(z-x)Q_j(z-y)\big| \,\big(|y-x|^\alpha+2^{-i\alpha}\big)\,dxdy   \\
&\lesssim \sum_{i;|i-j|\leq 4} 2^{-i\beta}\big(2^{-j\alpha} + 2^{-i\alpha}\big) \lesssim 2^{-j(\alpha+\beta)}.
\end{split} \end{equation*}
\end{Dem}

\medskip

Condition \eqref{EqWorkingHorse} is reminiscent of Gubinelli's notion of controlled path \cite{GubinelliControlled}. Recall from Proposition 35 in \cite{BB16} that for $f\in \mcC^{\alpha_1}$ and $g\in \mcC^{\alpha_2}$, with $\alpha_1,\alpha_2$ positive and $\alpha_1+\alpha_2\in (0,1)$, one has
$$
\Big|({\sf P}_fg)(y) - ({\sf P}_fg)(x) - f(x)\big(g(y)-g(x)\big)\Big| \lesssim |y-x|^{\alpha_1+\alpha_2}.
$$
It follows from Proposition \ref{PropWorkingHorse} above that ${\sf R}^\circ(f,g,h)\in\mcC^{\alpha_1+\alpha_2+\beta}$, for any $h\in\mcC^\beta$, with $\beta\in\bbR$. 

\medskip

Identity \eqref{EqGammayxIterateG} provides another example of a setting where Proposition \ref{PropWorkingHorse} applies, as it states that one has for any $\tau,\sigma\in\mcB^+$
\begin{equation*} \begin{split}
&{\sf g}_y\big(\tau/^+\sigma\big) - {\sf g}_x\big(\tau/^+\sigma\big)   \\
&= \sum_{n\geq 1} (-1)^n \hspace{-0.35cm}\sum_{\sigma<^+\sigma_n<^+\cdots<^+\tau} {\sf g}_x\big(\tau/^+\sigma_1\big)\cdots {\sf g}_x\big(\sigma_{n-1}/^+\sigma_n\big)\,\Big({\sf g}_y\big(\sigma_n/^+\sigma\big) - {\sf g}_x\big(\sigma_n/^+\sigma\big)\Big)   \\
&\quad+ {\sf g}_{yx}\big(\tau/^+\sigma\big),
\end{split} \end{equation*}   
with $\big|{\sf g}_{yx}\big(\tau/^+\sigma\big)\big| \lesssim |y-x|^{|\tau|-|\sigma|}$.   

\medskip

\begin{cor}   \label{CorWorkingHorse}
For any family $(h_\sigma)_{\sigma\in\mcB^+\backslash\{\bf 1\}}$, with $h_\sigma\in\mcC^{|\sigma|}$, the sum 
$$
\sum_{n\geq 1} (-1)^n \hspace{-0.35cm}\sum_{{\bf 1}<^+\sigma<^+\sigma_n<^+\cdots<^+\tau} {\sf R}^\circ\Big({\sf g}\big(\tau/^+\sigma_1\big)\cdots {\sf g}\big(\sigma_{n-1}/^+\sigma_n\big), {\sf g}(\sigma_n/\sigma), h_\sigma\Big)
$$
defines an element of $\mcC^{|\tau|}(\bbR^d)$.
\end{cor}

\bigskip

\subsection{Paracontrolled representations}

\begin{prop}
\label{PropDefnBracket}
Fix a regularity structure $\mathscr{T}$ and a model ${\sf M}=({\sf g}, {\sf \Pi})$ on $\mathscr{T}$. Define recursively on $|\tau|$ and $|\sigma|$ the families of real-valued functions $\big\{\Brac{\tau}^{\sf g}\big\}_{\tau\in\mcB^+}$ and $\big\{\Brac{\sigma}^{\sf M}\big\}_{\sigma\in\mcB}$ on $\bbR^d$, by the formulas
\begin{equation}
\label{EqDefnBrac1}
\begin{split}
{\sf g}(\tau)  &=  \sum_{{\bf 1}<^+\nu<^+\tau} {\sf P}_{{\sf g}(\tau/^+\nu)} \Brac{\nu}^{\sf g} + \Brac{\tau}^{\sf g}, \quad \tau\in\mcB^+,   \\
 {\sf \Pi}\sigma  &= \sum_{\mu<\sigma} {\sf P}_{{\sf g}(\sigma/\mu)}\Brac{\mu}^{\sf M} + \Brac{\sigma}^{\sf M}, \quad \sigma\in\mcB.
\end{split}
\end{equation}
Then $\Brac{\tau}^{\sf g}\in \mcC^{|\tau|}(\bbR^d)$, for all $\tau\in\mcB^+$, and $\Brac{\sigma}^{\sf M}\in \mcC^{|\sigma|}(\bbR^d)$, for all $\sigma\in\mcB$. Furthermore, the maps 
$$
{\sf M}\mapsto \Brac{\tau}^{\sf g}\in\mcC^{|\tau|}(\bbR^d), \quad {\sf M}\mapsto \Brac{\sigma}^{\sf M}\in\mcC^{|\sigma|}(\bbR^d),
$$
are continuous, for any $\tau\in\mcB^+$ and $\sigma\in\mcB$.
\end{prop}

\smallskip

\begin{Dem}
First we construct the family $\big\{\Brac{\tau}^{\sf g}; \tau\in\mcB^+\big\}$.   \vspace{0.1cm}

$\bullet$ The proof proceeds by induction on the homogeneity $|\tau|$ of $\tau$, starting with the case $\tau={\bf 1}$, for which we set $\Brac{{\bf 1}}^{\sf g} := {\sf g}({\bf 1}) = \textit{1}$, the constant function on $\bbR^d$, equal to $1$. Let $|\tau|>0$ and assume that the functions $\Brac{\sigma}^{\sf g}\in\Hoel^{|\sigma|}(\bbR^d)$ satisfying \eqref{EqDefnBrac1} have been constructed for any $\sigma\in\mcB^+$ with $|\sigma|<|\tau|$. Applying the two-parameter extension of the paraproduct operator ${\bfP}$ to identity \eqref{EqGammayxIterateG} with $\sigma = {\bf 1}$ and $<^+$ order, we have
$$
{\sf P}_{\textit{1}} {\sf g}(\tau) = \sum_{n=1}^\infty(-1)^{n-1}\sum_{\unit<^+\sigma_n<\cdots<\sigma_1<^+\tau}{\sf P}_{{\sf g}(\tau/^+\sigma_1)\cdots {\sf g}(\sigma_{n-1}/^+\sigma_n)} {\sf g}(\sigma_n) + {\bfP}_{x,y}\big({\sf g}_{yx}(\tau)\big).
$$
We used the fact that ${\sf P}_f \textit{1} = 0$, for any $f\in\Schw'(\Euc^d)$, to remove the zero-contribution from the $\sigma_n = {\bf 1}$ term in the sum. Note that 
$$
{\sf P}_{\textit{1}} {\sf g}(\tau) = {\sf g}(\tau) - \mathscr{S}{\sf g}(\tau),
$$ 
is the sum of ${\sf g}(\tau)$ and a smooth term depending continuously in any H\"older topology on ${\sf g}(\tau)\in L^\infty(\bbR^d)$. Expanding ${\sf g}(\sigma_n)$ by induction, we have
\begin{align*}
&\sum_{n=1}^\infty(-1)^{n-1}
\sum_{\unit<^+\sigma_n<^+\cdots<\sigma_1<^+\tau}
{\sf P}_{{\sf g}(\tau/^+\sigma_1)\cdots {\sf g}(\sigma_{n-1}/^+\sigma_n)}
{\sf g}(\sigma_n)\\
&= \sum_{n=1}^\infty(-1)^{n-1}
\sum_{\unit<^+\sigma_n<^+\cdots<^+\sigma_1<^+\tau}
{\sf P}_{{\sf g}(\tau/^+\sigma_1)\cdots {\sf g}(\sigma_{n-1}/^+\sigma_n)}
\Brac{\sigma_n}^{\sf g}   \\
&\quad+ \sum_{n=1}^\infty(-1)^{n-1}
\sum_{\unit<^+\sigma_{n+1}<^+\cdots<^+\sigma_1<^+\tau} 
{\sf P}_{{\sf g}(\tau/^+\sigma_1)\cdots {\sf g}(\sigma_{n-1}/^+\sigma_n)} 
\Big( {\sf P}_{{\sf g}(\sigma_n/^+\sigma_{n+1})}\Brac{\sigma_{n+1}}^{\sf g}\Big)   \\
&= \sum_{\unit<\sigma<\tau} {\sf P}_{{\sf g}(\tau/\sigma)} \Brac{\sigma}^{\sf g}   \\
&\quad+ \sum_{n=1}^\infty(-1)^{n-1} 
\sum_{\unit<^+\sigma_{n+1}<^+\cdots<^+\sigma_1<^+\tau} 
{\sf R}^\circ\Big({\sf g}(\tau/^+\sigma_1)\cdots {\sf g}(\sigma_{n-1}/^+\sigma_n), {\sf g}(\sigma_n/^+\sigma_{n+1}), \Brac{\sigma_{n+1}}^{\sf g}\Big).
\end{align*}
from a (wonderful) telescopic sum simplification. \textit{This is where something is happening}. Define then $\Brac{\tau}^{\sf g}$ by the formula
\begin{align*}
&\mathscr{S}{\sf g}(\tau) + {\bfP}_{x,y}\big({\sf g}_{yx}(\tau)\big)   \\
&+\sum_{n=1}^\infty(-1)^{n-1}\sum_{\unit<^+\sigma_{n+1}<^+\cdots<^+\sigma_1<^+\tau} {\sf R}^\circ\Big({\sf g}(\tau/^+\sigma_1)\cdots {\sf g}(\sigma_{n-1}/^+\sigma_n), {\sf g}(\sigma_n/^+\sigma_{n+1}), \Brac{\sigma_{n+1}}^{\sf g}\Big).
\end{align*}
It follows from the estimate $|{\sf g}_{yx}(\tau)|\lesssim|y-x|^{|\tau|}$, and Proposition \ref{PropParaprodEstimates} that ${\bfP}_{x,y}\big({\sf g}_{yx}(\tau)\big)\in\mcC^{|\tau|}(\bbR^d)$. Corollary \ref{CorWorkingHorse} takes care of the sum and shows that it defines an element of $C^{|\tau|}(\bbR^d)$.

\smallskip

$\bullet$ The proof of the regularity statement for $\Brac{\tau}^{\sf M}$, for $\tau\in\mcB$, proceeds by induction, similarly as above, using identity \eqref{EqPixIteratePi} giving ${\sf\Pi}^{\sf g}_x\tau$ in terms of ${\sf \Pi}$ only, as an input. Applying the two-parameter operator ${\bfP}$ to \eqref{EqPixIteratePi}, gives \eqref{EqDefnBrac1} for a choice of $\Brac{\tau}^{\sf M}$ equal to 
\begin{align*}
&\mathscr{S}{\sf \Pi}(\tau) + {\bfP}_{x,y}\big({\sf \Pi}_x^{\sf g}(\tau)(y)\big)   \\
&+\sum_{n=1}^\infty(-1)^{n-1}\sum_{\unit<\sigma_{n+1}<\cdots<\sigma_1<\tau} {\sf R}^\circ\Big({\sf g}(\tau/\sigma_1)\cdots {\sf g}(\sigma_{n-1}/\sigma_n), {\sf g}(\sigma_n/\sigma_{n+1}), \Brac{\sigma_{n+1}}^{\sf M}\Big).
\end{align*}
Since ${\sf\Pi}^{\sf g}_x\tau = {\sf\Pi}^{\sf g}_{x'}\widehat{{\sf g}_{x'x}}\tau = {\sf\Pi}^{\sf g}_{x'}\tau + \sum_{\sigma<\tau} {\sf g}_{x'x}(\tau/\sigma) {\sf\Pi}^{\sf g}_{x'}\sigma$, one can use Proposition \ref{PropReconstructionPC} to conclude that ${\bfP}_{x,y}\big({\sf \Pi}_x^{\sf g}(\tau)(y)\big)$ belongs to $\mcC^{|\tau|}(\bbR^d)$. 
\end{Dem}

\medskip

Recall from Example 2 in Section \ref{SectionBasicsRS} that given a model $\sf (g,\Pi)$ on $\mathscr{T}$, the concrete regularity structure $\mathscr{T}^+$ is endowed with an associated canonical model 
$$
{\sf M^g}:={\sf (g, \Pi^g) = (g,g)}.
$$ 
Remark that 
$$
\Brac{\cdot}^{\sf M^g} = \Brac{\cdot}^{\sf g},
$$ 
on $T^+$, so the above statement is really about $\Brac{\cdot}^{\sf M^g}$ and $\Brac{\cdot}^{\sf M}$. The proof makes it clear that the $\sf g$-brackets $\Brac{\tau}^{\sf g}$ depend only on $\sf g$. We extend by linearity the operators $\Brac{\cdot}^{\sf g}, \Brac{\cdot}^{\sf M}$ to $T^+$ and $T$, respectively. 

\medskip

\noindent \textbf{\textsf{Remark.}} {\it One can make the link with the setting introduced in \cite{BB16}, and give a different representation of the brackets under the assumption that we are given an operator $\bfI$ that acts on smooth functions and an abstract integration operator $\mcI : T\mapsto T$, on the regularity structure $\mathscr{T}$, together with a naive interpretation operator $\sf \Pi$, such that $\sf\Pi$ is multiplicative and ${\sf \Pi}(\mcI\tau)=\bfI({\sf \Pi}\tau)$, for all $\tau\in T$ -- see Section \ref{SectionSchauder}. Let then $\zeta : \bbR^d\mapsto\bbR$, stand for a smooth `noise' and $\circ$ stand for an element of $T$ such that ${\sf \Pi}(\circ)=\zeta$. We assign homogeneity $\alpha-\theta$ to $\circ$, and $|\tau|+\theta$ to any $\mcI\tau$, and $|\tau_1\cdots\tau_k|=|\tau_1|+\cdots+|\tau_k|$, for all $\tau_i\in T$. Set $g(\mcI\circ):={\sf \Pi}(\mcI\circ):=\bfI(\zeta)=:Z$. Denote by the unconventional sign $\circleddash$ the resonant operator from the paraproduct decomposition of a product -- see Appendix \ref{AppendixParapdcts}. Then
$$
{\sf \Pi}(\circ \mcI(\circ)) = \zeta Z = {\sf P}_Z\zeta + {\sf P}_\zeta Z + \circleddash(Z,\zeta),
$$
and we read on this expression that
$$
\Brac{\circ \mcI(\circ)}^{\sf M} = {\sf P}_\zeta Z + \circleddash(Z,\zeta).
$$
Compute $\Brac{\circ\mcI(\circ)^2}^{\sf M}$ as another example. We have 
$$
{\sf \Pi}\big(\circ\mcI(\circ)^2\big) = Z^2\zeta = {\sf P}_{Z^2}\zeta + 2{\sf P}_\zeta{\sf P}_ZZ + {\sf P}_\zeta\big(\circleddash(Z,Z)\big) + \circleddash\big(2{\sf P}_ZZ+\circleddash(Z,Z),\zeta\big).
$$
To make the link with the defining relation \eqref{EqDefnBrac1} for $\Brac{\circ\mcI(\circ)^2}^{\sf M}$, we use the corrector operator ${\sf C}$ and the operator $\sf S$ from \cite{BB16}. This gives for ${\sf \Pi}\big(\circ\mcI(\circ)^2\big)$ the expression
\begin{equation*}
\begin{split}
 &{\sf P}_{Z^2}\zeta + 2{\sf P}_Z\big({\sf P}_\zeta Z\big) + 2{\sf S}(\zeta,Z,Z) + {\sf P}_\zeta\big(\circleddash(Z,Z)\big) + 2Z\circleddash(Z,\zeta) + 2{\sf C}(Z,Z,\zeta) + \circleddash\big(\circleddash(Z,Z),\zeta\big)   \\
&= {\sf P}_{Z^2}\zeta + {\sf P}_Z\big(2{\sf P}_\zeta Z + \circleddash(Z,\zeta)\big) + \Big\{2{\sf S}(\zeta,Z,Z) + {\sf P}_\zeta\big(\circleddash(Z,Z)\big) + 2{\sf P}_{\circleddash(Z,\zeta)}Z + 2\circleddash\big(Z,\circleddash(Z,\zeta)\big)   \\
&\hspace{5.3cm}+2{\sf C}(Z,Z,\zeta) + \circleddash\big(\circleddash(Z,Z),\zeta\big)\Big\},
\end{split}
\end{equation*}
so the term inside the brackets $\{\cdots\}$ defines $\Brac{\circ\mcI(\circ)^2}^{\sf M}$. As can be seen from these examples, these expressions of the brackets using the operators from \cite{BB16} quickly get seemingly complicated.}

\medskip

\begin{prop}
We have
$$
{\sf g}(\tau/\sigma) = \sum_{\eta\in\mcB,\,\sigma<\eta<\tau} {\sf P}_{{\sf g}(\tau/\eta)}\Brac{\eta/\sigma}^{\sf g} + \Brac{\tau/\sigma}^{\sf g},
$$
for all $\sigma, \tau\in\mcB$ with $\sigma<\tau$.
\end{prop}

\medskip

\begin{Dem}
With $\tau/\sigma = \sum_{\theta\in\mcB^+} (\Delta\tau)^{\sigma\theta}\,\theta $, and $\theta/^+\rho = \sum_{\kappa\in\mcB^+}(\Delta^+\theta)^{\rho\kappa}\kappa$, we have on the one hand, from Proposition \ref{PropDefnBracket},
\begin{equation*}
\begin{split}
{\sf g}(\tau/\sigma)  &=  \sum_{\theta\in\mcB^+} (\Delta\tau)^{\sigma\theta} {\sf g}(\theta) =  \sum_{\theta,\rho\in\mcB^+,\unit<^+\rho<^+\theta} (\Delta\tau)^{\sigma\theta} \, {\sf P}_{{\sf g}(\theta/^+\rho)} \Brac{\rho}^{\sf g} + \sum_{\theta\in\mcB^+} (\Delta\tau)^{\sigma\theta} \, \Brac{\theta}^{\sf g}   \\
		 						   &=  \sum_{\theta,\rho,\kappa\in\mcB^+,\,\unit<^+\rho<^+\theta} (\Delta\tau)^{\sigma\theta} (\Delta^+\theta)^{\rho\kappa} \, {\sf P}_{{\sf g}(\kappa)} \Brac{\rho}^{\sf g} + \Brac{\tau/\sigma}^{\sf g},
\end{split}
\end{equation*}
and on the other hand, and since ${\sf P}_f {\it 1} = 0$, for any $f\in\mcS'(\bbR^d)$,
$$
\sum_{\sigma<\eta<\tau} {\sf P}_{{\sf g}(\tau/\eta)} \Brac{\eta/\sigma}^{\sf g}  =  \sum_{\sigma<\eta<\tau,\kappa,\rho\in\mcB^+,\mathbf{1}<^+\rho} (\Delta\tau)^{\eta\kappa} (\Delta\eta)^{\sigma\rho} \, {\sf P}_{{\sf g}(\kappa)} \Brac{\rho}^{\sf g}.
$$
The statement then follows from the right comodule identity \eqref{EqRightComodule}, such as expressed in coordinates in identity \eqref{EqRightComoduleComponents}, and the structural assumption \eqref{EqDelta} on the splitting map $\Delta$.
\end{Dem}

\medskip

\begin{thm}
\label{ThmReconstructionPC}
Fix a regularity exponent $\gamma\in\bbR$, and a model ${\sf M}={\sf (g,\Pi)}$ on the regularity structure $\mathscr{T}$. One can associate to any modelled distribution 
$$
\bsf=\sum_{\tau\in\mcB; |\tau|<\gamma} f^\tau \tau \in \mcD^\gamma(\mathscr{T}, {\sf g}),
$$ 
a distribution $\Brac{\bsf}^{\sf M} \in \mcC^\gamma(\bbR^d)$ such that one defines a reconstruction $\textsf{\textbf{R}}\bsf$ of $\bsf$ setting
\begin{align}\label{2: eq reconst to PD}
\textsf{\textbf{R}}\bsf  :=  \sum_{\tau\in\mcB; |\tau|<\gamma} {\sf P}_{f^\tau} \Brac{\tau}^{\sf M} + \Brac{\bsf}^{\sf M}.
\end{align}
Each coefficient $f^\tau$ also has a representation
\begin{align}\label{2: eq coefficients to PD}
f^\tau  =  \sum_{\tau<\mu; |\mu|<\gamma} {\sf P}_{f^\mu} \Brac{\mu/\tau}^{\sf g} + \Brac{f^\tau}^{\sf g},
\end{align}
for some $\Brac{f^\tau}^{\sf g}\in\mcC^{\gamma-|\tau|}(\bbR^d)$. Moreover, the map 
$$
\bsf \mapsto \Big(\Brac{\bsf}^{\sf M}, \big(\Brac{f^\tau}^{\sf g}\big)_{\tau\in\mcB} \Big)
$$
from $\mcD^\gamma(\mathscr{T}, {\sf g})$ to $\mcC^\gamma(\bbR^d)\times \prod_{\tau\in\mcB} \mcC^{\gamma-|\tau|}(\bbR^d)$, is continuous.
\end{thm}

\medskip

\begin{Dem}
Recall from Proposition \ref{PropReconstructionPC} that, there exists a function $g\in\mcC^\gamma(\bbR^d)$ such that
\begin{equation}
\label{EqDecompositionReconstruction}
{\bfP}_{x,y}\Big(\big({\sf\Pi}_x^{\sf g}\bsf(x)\big)(y)\Big) - g\mathbf{1}_{\gamma>0}
\end{equation}
is the reconstruction $\textsf{\textbf{R}}\bsf$. We have from identity \eqref{EqPixIteratePi} giving ${\sf\Pi}^{\sf g}_x$ in terms of $\sf g$ and $\sf \Pi$, the finite expansion
$$
\big({\sf\Pi}^{\sf g}_x\bsf(x)\big)(\cdot) = \sum_{n=0}^\infty (-1)^n \sum_{\sigma_n<\cdots<\sigma_1<\sigma_0}\,\Big(f^{\sigma_0}{\sf g}\big(\sigma_0/\sigma_1\big)\cdots {\sf g}\big(\sigma_{n-1}/\sigma_n\big)\Big)(x)\,({\sf \Pi}\sigma_n)(\cdot).
$$
So applying the two-parameter paraproduct operator $\bfP$ on both sides, and using the same (fantastic) telescopic sum as in the proof of Proposition \ref{PropDefnBracket}, we get, with an obvious notation,
\begin{equation*}
\begin{split}
&{\bfP}_{x,y}\Big(\big({\sf\Pi}^{\sf g}_x\bsf(x)\big)(y)\Big)  = (\mcC^\infty)+  \sum_{\sigma_0\in\mcB} {\sf P}_{f^{\sigma_0}} \Brac{\sigma_0}^{\sf M}   \\
&+ \sum_{n=1}^\infty (-1)^n \sum_{\sigma_{n+1}<\sigma_n<\cdots<\sigma_1<\sigma_0} \, {\sf R}^\circ\Big(f^{\sigma_0}{\sf g}\big(\sigma_0/\sigma_1\big)\cdots {\sf g}\big(\sigma_{n-1}/\sigma_n\big), {\sf g}(\sigma_n/\sigma_{n+1}), \Brac{\sigma_{n+1}}^{\sf M}\Big).
\end{split}
\end{equation*}
From \eqref{EqIncrementMdF}, for each $\sigma\in\mcB$ we have
\begin{align*}
&f^\sigma(y)-f^\sigma(x)\\
&=\big(\bsf(y)-\widehat{{\sf g}_{yx}}\bsf(x)\big)^\sigma+\sum_{\sigma<\sigma_0}f^{\sigma_0}(x){\sf g}_{yx}(\sigma_0/\sigma)\\
&=\big(\bsf(y)-\widehat{{\sf g}_{yx}}\bsf(x)\big)^\sigma\\
&\quad+\sum_{n\ge0}(-1)^n\sum_{\sigma<\sigma_n<\cdots<\sigma_1<\sigma_0}\big(f^{\sigma_0}{\sf g}(\sigma_0/\sigma_1)\cdots{\sf g}(\sigma_{n-1}/\sigma_n)\big)(x)\big({\sf g}_y(\sigma_n/\sigma)-{\sf g}_x(\sigma_n/\sigma)\big).
\end{align*}
Proposition \ref{PropWorkingHorse} applies and tells us that the sum of the ${\sf R}^\circ$-terms defines an element of H\"older regularity $(\gamma-|\sigma_{n+1}|)+|\sigma_{n+1}|=\gamma$. The claim of the theorem on $\textsf{\textbf{R}}\bsf$ comes from this fact and identity \eqref{EqDecompositionReconstruction}. To get the paracontrolled representation of $f^\tau$, note from Lemma \ref{LemmaExModelledFunction} that $\bsf/\tau = \sum_{\mu\ge\tau; |\mu|<\gamma} f^\mu (\mu/\tau)\in \mcD^{\gamma-|\tau|}(\mathscr{T}^+,{\sf g})$, and apply the result just proved to the reconstructions of the modelled distribution.
\end{Dem}

\medskip

Theorem \ref{ThmReconstructionPC} refines over the paraproduct-based construction of the reconstruction operator given by Gubinelli, Imkeller and Perkowski in \cite{GIP15}, where $\textbf{\textsf{R}}\bsf$ is proved to be of the form ${\bfP}_{x,y}\big(({\sf\Pi}^{\sf g}_x\bsf(x))(y)\big)$, up to a $\mcC^\gamma(\bbR^d)$ term. See our extension in Proposition \ref{PropReconstructionPC} above. The point of our refined representation of the family $\big(\textbf{\textsf{R}}\bsf, f^\tau\big)$ as a paracontrolled system lies in Theorem \ref{ThmMain2}, proved in the next section. It parametrizes the class of ``admissible" models used for the study of singular stochastic PDEs, in terms of the brackets $\Brac{\tau}^{\sf M}$, with $|\tau|\leq 0$. The forthcoming work \cite{BH2} will give a similar description of more general models $\sf (g,\Pi)$ and modelled distributions in $\mcD^\gamma(\mathscr{T},{\sf g})$, in terms only of the bracket data, extending the main result of \cite{MP18} on Besov-type characterizations $\mcD^\gamma(\mathscr{T},{\sf g})$.

\ssk

One advantage of the explicit construction of the reconstruction operator given by Theorem \ref{ThmReconstructionPC} is that this representation is flexible enough to work in other functional settings than the present $B^\gamma_{\infty\infty}$-type space $\mcD^\gamma(\mathscr{T}, {\sf g})$. The continuity properties of the paraproduct operator on Besov, Triebel-Lizorkin or Sobolev-Slobodeckij spaces are well-known, and allow for a direct approach to reconstruction in these spaces, in the line of the recent works \cite{HL16, HS17, LPT16, ST18}.

\medskip

Before giving the next statement, note that the restriction to $T_{\leq 0}$ of the splitting $\Delta$ turns 
$$
\mathscr{T}_{\leq 0} := \Big((T^+,\Delta^+), (T_{\leq 0},\Delta)\Big)
$$ 
into a regularity structure. The next statement is essentially contained in Proposition 3.31 from \cite{Hai14}; we give the details here to provide a self-contained document.

\medskip

\begin{cor}
\label{CorUniqueExtension}
Assume we are given a map ${\sf g}:\bbR^d\to G^+$ such that the bound \eqref{EqEstimateGammayx} is satisfied. Let a family $\big(\Brac{\tau}\in\mcC^{|\tau|}(\bbR^d)\big)_{\tau\in\mcB, |\tau|\leq 0}$ be given. For any $\tau\in\mcB$ with $|\tau|\leq 0$, set 
$$
{\sf \Pi}\tau := \sum_{\sigma<\tau} {\sf P}_{{\sf g}(\tau/\sigma)}\Brac{\sigma} + \Brac{\tau}.
$$
Then $({\sf g},\sf \Pi)$ is a model on the regularity structure $\mathscr{T}_{\leq 0}$, and it has a unique extension into a model on $\mathscr{T}$.
\end{cor}

\medskip

\begin{Dem}
$\bullet$ Pick a basis vector $\tau\in\mcB$ with $|\tau|\leq 0$, and assume that $\sf(g,\Pi)$ is a model on the sector $T_{<|\tau|}$. Set for all $x\in\bbR^d$
$$
\bsh_\tau(x) := \sum_{\sigma<\tau}\, {\sf g}_x(\tau/\sigma)\sigma;
$$
this defines a modelled distribution in $\mcD^{|\tau|}({\mathscr{T}},{\sf g})$. Then the bound $|\langle{\sf \Pi}_x^{\sf g}\tau,\varphi_x^\lambda\rangle|\lesssim\lambda^{|\tau|}$ is equivalent to that ${\sf \Pi}\tau$ is (one of) the reconstructions of $\bsh_\tau$. From Theorem \ref{ThmReconstructionPC}, the distribution 
$$
\sum_{\sigma<\tau} {\sf P}_{{\sf g}(\tau/\sigma)}\Brac{\sigma} + \Brac{\bsh_\tau}
$$
is a reconstruction of $\bsh_\tau$. Since 
$$
{\sf \Pi}\tau - \bsh_\tau = \Brac{\tau} - \Brac{\bsh_\tau} \in \mcC^{|\tau|}(\bbR^d),
$$
the distribution ${\sf \Pi}\tau$ appears then as another reconstruction of $\bsh_\tau$.

\medskip

$\bullet$ If one picks now a basis vector $\mu\in\mcB$, with $|\mu|>0$, then $\bsh_\mu\in\mcD^{|\mu|}(\mathscr{T},{\sf g})$ has a unique reconstruction, equal to ${\sf \Pi}\mu$, that is characterized by the data $\big({\sf\Pi}_x^{\sf g}\sigma, {\sf g}_x(\mu/\sigma) ; x\in\bbR^d, \sigma<\mu\big)$, from the defining property \eqref{EqReconstructionCondition} of a reconstruction. An elementary induction then shows the existence of a unique extension of $\sf \Pi$ to $T$ that satisfies the property ${\sf \Pi}\tau = \textsf{\textbf{R}}\bsh_\tau$, for every $\tau\in\mcB$ with positive homogeneity.
\end{Dem}

\bigskip

\section{Parametrization of the set of admissible models}
\label{SectionParametrization}

\subsection{Usual models}
\label{SectionUsual}

We introduce in this section a notion of usual model on a concrete regularity structure, motivated by some identity satisfied by ${\sf g}(\tau/\bsX^k)$ in the usual setting; see Equation \eqref{usualness of model} below. Its introduction is motivated by the fact that usual models $\sf (g,\Pi)$ are entirely determined by the $\sf \Pi$ map, under a mild structure assumption on $T^+$ and $\Delta$. The definition of a usual model requires that we work with concrete regularity structures where $T$ and $T^+$ are related with one another, unlike the results of the previous section.

\medskip

Let $\mathscr{T} = \big((T^+,\Delta^+), (T,\Delta)\big)$ be a concrete regularity structure. If $T$ contains the \emph{usual} polynomial structure $T(X)$,  one can expand the coproduct $\Delta\tau$ of any $\tau\in\mathcal{B}\setminus\{\bsX^k\}_k$, as 
$$
\Delta\tau = \mathring{\Delta}\tau + \sum_{k\in\mathbb{N}^d}\frac{\bsX^k}{k!}\otimes \bfD^k\tau,
\quad
\mathring{\Delta}\tau := \sum_{\sigma\neq \bsX^k}\sigma\otimes(\tau/\sigma),
$$
where 
$$
\bfD^k\tau := k!(\tau/\bsX^k).
$$
Applying ${\sf \Pi}\otimes{\sf g}_x^{-1}$, we have
$$
{\sf\Pi}_x^{\sf g}\tau = \mathring{{\sf\Pi}}_x^{\sf g}\tau + \sum_{k\in\mathbb{N}^d}\frac{(\cdot)^k}{k!}\, {\sf g}_x^{-1}(\bfD^k\tau),
\quad
\mathring{\sf\Pi}_x^{\sf g}:=({\sf \Pi}\otimes{\sf g}_x^{-1})\mathring{\Delta}.
$$
Setting
$$
{\sf f}_x(\bfD^k\tau):=-\sum_{\ell\in\mathbb{N}^d}\frac{x^\ell}{\ell!}\,{\sf g}_x^{-1}(\bfD^{k+\ell}\tau),
$$
or equivalently, 
$$
{\sf g}_x^{-1}(\bfD^k\tau):=-\sum_{\ell\in\mathbb{N}^d}\frac{(-x)^\ell}{\ell!}\,{\sf f}_x(\bfD^{k+\ell}\tau),
$$ 
gives a Taylor-like expansion formula for ${\sf\Pi}^{\sf g}_x\tau$, under the form of the identity
$$
{\sf\Pi}^{\sf g}_x\tau = \mathring{\sf\Pi}_x^{\sf g}\tau - \sum_{k\in\mathbb{N}^d}\frac{(\cdot-x)^k}{k!}\,{\sf f}_x(\bfD^k\tau).
$$
Since the derivatives $\partial_y^k({\sf\Pi}^{\sf g}_x\tau)(y)$ vanishes at $y=x$ for any $|k|<|\tau|$, one has
\begin{equation}\label{eq f_x(tau/X^k)}
{\sf f}_x(\bfD^k\tau)=\mathbf{1}_{|k|<|\tau|}\partial_y^k(\mathring{\sf\Pi}_x^{\sf g}\tau)(y)\big|_{y=x}.
\end{equation}
Given $\alpha\in\mathbb{R}$, define a linear projection map $P_{>\alpha} : T\mapsto T$ setting
$$
P_{>\alpha}(\tau) := \tau \, {\bf 1}_{|\tau|>\alpha},
$$ 
for every $\tau\in\mcB$.

\medskip

\begin{lem}
For any $\tau\in\mathcal{B}\setminus\{\bsX^k\}_k$, one has
\begin{equation}\label{eq g(tau/X^k)}
{\sf g}_x(\bfD^k\tau) = \Big(\partial_y^k\big\{(\mathring{\sf\Pi}_x^{\sf g} P_{>|k|}\otimes{\sf g}_x)\Delta \tau\big\}(y)\Big) \Big|_{y=x}.
\end{equation}
\end{lem}

\medskip

\begin{Dem}
It suffices to show that
\begin{equation}\label{eq tau/X^k}
{\sf g}_x(\bfD^k\tau)=\sum_{\sigma\neq \bsX^\ell} {\sf f}_x(\bfD^k\sigma){\sf g}_x(\tau/\sigma);
\end{equation}
we get \eqref{eq g(tau/X^k)} by inserting \eqref{eq f_x(tau/X^k)} into \eqref{eq tau/X^k}. We start from the formula
$$
\Delta^+(\tau/\bsX^{k+\ell}) = \sum_{\sigma\neq \bsX^m}(\sigma/\bsX^{k+\ell})\otimes(\tau/\sigma) + \sum_m
\begin{pmatrix}
k+\ell+m\\
m
\end{pmatrix}
\bsX^m\otimes(\tau/\bsX^{k+\ell+m}).
$$
Since $\tau/\bsX^{k+\ell}\in T^+\backslash \langle\unit\rangle$, applying ${\sf g}_x^{-1}\otimes{\sf g}_x$ to the preceding identity gives
$$
0 = \sum_{\sigma\neq \bsX^m}{\sf g}_x^{-1}(\bfD^{k+\ell}\sigma){\sf g}_x(\tau/\sigma) + \sum_m\frac{(-x)^m}{m!}{\sf g}_x(\bfD^{k+\ell+m}\tau),
$$
that is
$$
0=
-\sum_m\frac{(-x)^m}{m!}
\left(\sum_{\sigma\neq \bsX^m}{\sf f}_x(\bfD^{k+\ell+m}\sigma){\sf g}_x(\tau/\sigma)-{\sf g}_x(\bfD^{k+\ell+m}\tau)\right).
$$
Identity \eqref{eq tau/X^k} is obtained as a consequence, since
\begin{align*}
0
&=
\sum_{\ell,m}\frac{x^\ell}{\ell!}\frac{(-x)^m}{m!}
\left(\sum_{\sigma\neq \bsX^m}{\sf f}_x(\bfD^{k+\ell+m}\sigma){\sf g}_x(\tau/\sigma)-{\sf g}_x(\bfD^{k+\ell+m}\tau)\right)\\
&=
\sum_{\sigma\neq \bsX^m}{\sf f}_x(\bfD^k\sigma){\sf g}_x(\tau/\sigma)-{\sf g}_x(\bfD^k\tau).
\end{align*}
\end{Dem}

\medskip

We use in the present work the \emph{bounded} polynomial structure rather than the usual polynomial structure. We work with concrete regularity structures for which the following assumptions hold true.

\medskip

\noindent \textsf{\textbf{Assumption (A)}}   \hspace{0.1cm}
{\it
The bounded polynomial structure $T(X)=\langle \bsX_e^k \rangle_{e,k}$ is contained in $T$, and the polynomial ring $T^+(X)=\langle \bsX_{e_1}^{k_1}\cdots\bsX_{e_n}^{k_n} \rangle_{e_1,\dots,e_n,k_1,\dots,k_n}$ is included in $T^+$.
}

\medskip

We do not make a difference in the notations between the two copies in $T$ and $T^+$ of the bounded polynomial structure.

\medskip

\begin{defn}
Let $\mathscr{T}$ be a concrete regularity structure satisfying Assumption \textsf{\textbf{(A1)}}. We say that the \textsf{\textbf{model}} $({\sf g}, {\sf \Pi})$ is \textsf{\textbf{usual}} if one has ${\sf g}_x(\bsX_e^k)=({\sf \Pi}\bsX_e^k)(x)=x_e^k$, and
\begin{equation}\label{usualness of model}
{\sf g}_x(\bfD_e^k\tau) = \Big(\partial_y^k\left\{\phi_e\big(\mathring{\sf\Pi}_x^{\sf g} P_{>|k|}\otimes{\sf g}_x\big)\Delta \tau\right\} (y)\Big)\Big|_{y=x}.
\end{equation}
for any $\tau\in\mathcal{B}\setminus\{\bsX_e^k\}_{e,k}$, where $\bfD_e^k\tau:=k!(\tau/\bsX_e^k)$.
\end{defn}

\bigskip

\subsection{Abstract integration operator and admissible models}
\label{SectionSchauder}

Fix a positive regularity exponent $\theta$, and let $\mathscr{T}$ be a concrete regularity structure. Assume for simplicity that 
$$
\beta_0 > -\theta,
$$ 
so all the elements of $T$ have homogeneity strictly greater than $-\theta$. We consider in this section concrete regularity structures $\mathscr{T}$ equipped with an abstract integration operator $\mcI$, that is a regularity structure counterpart of an operator $\bfI$ that is typically an integral operator given by a kernel that is singular on the diagonal, such as the Green function of a differential operator. The exponent $\theta$ quantifies the regularizing properties of the operator $\bf I$ in the H\"older or Besov scale.

\smallskip

\begin{Remark*}
The dynamical $\Phi_3^4$ equation
$$
\partial_t\Phi=\Delta\Phi-\Phi^3+\xi
$$
seems not to satisfy the above assumption. Indeed, one should choose $\beta_0=-\frac52-\epsilon$, and $\theta=2$ for the heat kernel in any dimension. However, if we decompose $\Phi=\Psi+v$, where $\partial_t\Psi=\Delta\Psi+\xi$ and
$$
\partial_tv=\Delta v-(v+\Psi)^3,
$$
then one can choose $\beta_0=3(-\frac12-\epsilon)$ instead, so the equation for $v$ satisfies $\beta_0>-\theta$. A general da Prato-Debussche trick is described in Section 6 of \cite{BCCH18}, that allows to set the study of a generic subcritical singular partial differential equation, within the setting of regularity structures, under the assumption $\beta_0>-\theta$.
\end{Remark*}

\bigskip

\hfill \textcolor{gray}{\textsf{\textbf{Integration operator}}}   

\bigskip

Let $K_n : \bbR^d\times\bbR^d \mapsto\bbR$, be a sequence of kernels on $\bbR^d$, with support in 
$$
\big\{(x,y)\in\bbR^d\times\bbR^d; |y-x|\leq 2^{-n}\big\},
$$ 
and such that one has, for all $n\in\bbN$ and $x,y\in\bbR^d$,
\begin{equation}
\label{EqConditionKn}
\big|\partial_x^k\partial_y^\ell K_n(x,y)\big|  \leq  C_{k,\ell} \, 2^{-n(\theta+\epsilon-d-|k|-|\ell|)},
\end{equation}
for some (small) positive $\epsilon$. (This $\epsilon$ is only needed in the proof of Lemma \ref{LemConstructionAdmissibleModels}; see the remark following that lemma.)
The converging sum
$$
K(x,y) = \sum_{n\geq 0} K_n(x,y)
$$
defines a kernel 
$$
K : \bbR^d\times\bbR^d\backslash\{(x,x);x\in\bbR^d\} \mapsto\bbR,
$$ 
and, for each $x\in\bbR^d$, an integration map
$$
(\bfI\varphi)(x) := \int_{\bbR^d} K(x,y)\varphi(y)\,dy,
$$
for $\varphi\in\mcD(\bbR^d\backslash\{x\})$. The archetypal example is the smoothly localized Green kernel 
$$
K(x,y) = \chi(|y-x|)\,|y-x|^{2-d},
$$ 
for $\chi$ a smooth real-valued function with compact support identically equal to $1$ in a neighbourhood of $0$, in dimension at least $d\geq 3$, for which one can take any $\theta<2$. The associated integration map sends any $\mcC^\beta(\bbR^d)$, into $\mcC^{\beta+2}(\bbR^d)$, for $\beta\notin\bbZ$ -- these are Schauder estimates. Note however that $(\bfI\zeta)(x)$ is not defined for a generic distribution $\zeta$.

\medskip

\begin{lem}
\label{LemmeDefnIDistribution}
Let $\{\zeta_x\}_{x\in\bbR^d}\subset\mcS'(\bbR^d)$ be a family of distributions. If there exist $\alpha\in\Euc$ and a positive constant $C$ such that one has
$$
\big|\langle\zeta_x,\varphi_x^\lambda\rangle\big|  \leq  C\delta^\lambda,
$$
uniformly over $\varphi\in\mcC^r(\bbR^d)$, with unit norm and support in the unit ball, $\lambda\in(0,1]$ and $x\in\Euc^d$, then the sum
\begin{equation}
\label{EqDefnIDistribution}
(\bfI_k^e(\zeta_x))(x) := \Big\langle \zeta_x, \partial_x^k\big(\phi_e(x) K(x,\cdot)\big) \big\rangle := \sum_{n\geq 0} \Big\langle \zeta_x, \partial_x^k\big(\phi_e(x) K_n(x,\cdot)\big)\Big\rangle
\end{equation}
converges for any $|k|<\alpha+\theta$, $e\in E$ and $x\in\bbR^d$.
\end{lem}

\medskip

\begin{Dem}
Pick $x\in\bbR^d$. Let $\varphi$ be any smooth function with support in $\{y\in\bbR^d;|y-x|<\lambda\}$, such that one has
$$
\sup_{|\ell|\le r}\lambda^{d+|\ell|}\|\partial^\ell\varphi\|_{L^\infty}  \leq  1,
$$
for some $\lambda\in(0,1]$. Since $\varphi_{x,\lambda}(y) := \lambda^d\varphi(x+\lambda y)$ has unit norm in $\mcC^r(\bbR^d)$ and $\varphi=(\varphi_{x,\lambda})_x^\lambda$, we have $|\langle\zeta_x,\varphi\rangle|\le C\lambda^\alpha$, from the assumption of the lemma.
Pick $k\in\bbN^d$. Since $\varphi(y) := \partial_x^k(\phi_e(x)K_n(x,y))$ is supported in $\{y\in\bbR^d;|y-x|<2^{-n}\}$ and $\|\partial^\ell\varphi\|_{L^\infty}\lesssim2^{(d+|k|+|\ell|-\theta)n}$, we thus have 
\begin{align}\label{Bound zeta_x(K_n)}
\Big|\Big\langle \zeta_x, \partial_x^k\big(\phi_e(x) K_n(x,\cdot)\big)\Big\rangle\Big| \lesssim 2^{(|k|-\alpha-\theta)n},
\end{align}
and a converging sum in \eqref{EqDefnIDistribution} if $|k|<\alpha+\theta$.
\end{Dem}

\smallskip

Note that we cannot even make sense of $\int_{\bbR^d}K(z,y)\zeta_x(y)\,dy$, for $z\neq x$.
Were we able to define that function as a regular function of $z$, it would have a regularity structure lift in the canonical polynomial structure. Lemma \ref{LemmeDefnIDistribution} allows to define an avatar for the lift at point $x$ only of the non-existing function $\big((\phi_e\bfI)\zeta_x\big)(\cdot)$, under the form of the quantity
$$
\sum_{e\in E, |k|<\alpha+\theta}\, \frac{\bsX_e^k}{k!}\, \big((\phi_e\bfI)\zeta_x\big)(x).
$$
It follows from Lemma \ref{LemmeDefnIDistribution} and the assumption $\beta_0>-\theta$, that one can make sense of $\bfI({\sf \Pi}\tau)(x)$ for any $x\in\bbR^d$, under the form of the converging sum
$$
\bfI({\sf \Pi}\tau)(x) := \sum_{n\geq 0} \big\langle {\sf\Pi}\tau,K_n(x,\cdot) \big\rangle.
$$

\bigskip

\hfill \textcolor{gray}{\textsf{\textbf{Regularity structures with an abstract integration operator}}}   \vspace{0.1cm}

\medskip

In addition to Assumption \textbf{\textsf{(A)}}, we make the following set of assumptions on the concrete regularity structure $\mathscr{T}$.

\medskip

\noindent \textsf{\textbf{Assumption (B)}}   \vspace{0.15cm}
\hspace{0.1cm}
{\it The sets $T^+$ and $T$ are related via the integral operators in the following sense.
\begin{itemize}
   \item There exist operators $\mcI^{e+}_k : T\mapsto T^+$, indexed by $e\in E$ and $k\in\bbN^d$, with positive homogeneities  
   $$
   |X_e^i|= 1,  \quad \big| \mcI^{e+}_k\tau\big| = |\tau|+\theta-|k|.
   $$   \vspace{0.1cm}   
   
   \item One has
   $$
   \copro\unit=\unit\otimes\unit, \quad \copro X_e^i = X_e^i\otimes\unit + \unit\otimes X_e^i,
   $$ 
   and the operators $\Delta$ and $\Delta^+$ are related by the \emph{intertwining relations}  
   \begin{equation}
   \label{EqIntertwiningDeltaDelta+}
   \copro(\mcI_k^{e+}\tau) = (\mcI_k^{e+}\otimes\iden)\comul\tau + \sum_{\ell\in\bbN^d} \frac{\bsX_e^\ell}{\ell!}\otimes \mcI_{k+\ell}^{e+} \tau.
\end{equation} \vspace{0.1cm}
\end{itemize}

In addition, $T$ satisfies the following assumptions.
\begin{itemize}
   \item There exists an operator $\mcI: T\mapsto T$, with 
   $$
   \big|\mcI\tau\big| = |\tau|+\theta.
   $$ 
\vspace{0.1cm}   
   \item For any $\tau\in \mcB$, one has
   \begin{equation}
   \label{EqDefnDeltaIkTau}
   \Delta(\mcI\tau) = (\mcI\otimes\iden)\Delta\tau + \sum_{e\in E,\ell\in\bbN^d} \frac{\bsX_e^\ell}{\ell!}\otimes \mcI_{\ell}^{e+}\tau.
   \end{equation}
\end{itemize}}   

\medskip

Note that identity \eqref{EqDefnDeltaIkTau} identifies $\mcI^{e+}_{k}\tau$ as $\mcI\tau/\bsX_e^k$, for any $k\in\bbN^d, e\in E$. The operators $\mcI^{e+}_k$ are the regularity structure counterparts of the operators $\partial^k(\phi_e\bfI)$. Note that the restrictions on the index sets in identities \eqref{EqIntertwiningDeltaDelta+} and \eqref{EqDefnDeltaIkTau} to indices $\ell$ with $ |k|+|\ell|<|\tau|+\theta$, are redundant with the fact that $\mcI^{e+}_{k+\ell}$ is null on $T_\beta$, for $\beta\le-|k|-|\ell|$. In applications to the study of stochastic PDEs with derivatives of unknown functions, such as the KPZ equation, we can also assume the existence of other operators $\mcI_k:T\to T$, associated with the integration operator $\partial^k\bfI$.

\medskip

\begin{prop}
\label{PropRealizationPiPi+}
Let $\sf (g,\Pi)$ be a usual model on $\mathscr{T}$.
We assume the commutation rule
\begin{align}\label{commutation I and Pi}
{\sf \Pi}(\mcI\tau) = \bfI({\sf \Pi}\tau).
\end{align}
Then, the usual property \eqref{usualness of model} holds for any $\tau\in\mcI\mcB$ if and only if, for every $\tau\in \mcB$, and $x\in\bbR^d$, one has
\begin{align}
\label{EqDefnAdmissibleGBis}
{\sf g}_x(\mcI^{e+}_k\tau) &=\sum_{\sigma\le\tau; |k|<|\sigma|+\theta}\, {\sf g}_x(\tau/\sigma)\, \bfI_k^e\big({\sf\Pi}^{\sf g}_x\sigma\big)(x)
\end{align}
\end{prop}

\medskip

\begin{Dem}
Since $\mcI_k^{e+}\tau=\bfD_e^k\mcI\tau$,
\begin{align*}
{\sf g}_x(\mcI_k^{e+}\tau) 
&= \Big(\partial_y^k\left\{\phi_e\big(\mathring{\sf\Pi}_x^{\sf g} P_{>|k|}\otimes{\sf g}_x\big)\Delta \mcI\tau\right\} (y)\Big)\Big|_{y=x}.  \\
&=\partial_y^k\Big\{\phi_e\Big(
\sum_{\sigma;|\sigma|+\theta>|k|}(\mathring{\sf\Pi}_x^{\sf g} \mcI\sigma){\sf g}_x(\tau/\sigma)
\Big)\Big\}(y)\Big|_{y=x}   \\
&=\partial_y^k\Big\{\phi_e\Big(
\sum_{\sigma,\eta;|\sigma|+\theta>|k|}({\sf\Pi}\mcI\eta){\sf g}_x^{-1}(\sigma/\eta){\sf g}_x(\tau/\sigma)
\Big)\Big\}(y)\Big|_{y=x}   \\
&=\bfI_k^e\left(\sum_{\sigma,\eta;|\sigma|>|k|-\theta}({\sf\Pi}\eta){\sf g}_x^{-1}(\sigma/\eta){\sf g}_x(\tau/\sigma) \right)(x)   \\
&=\bfI_k^e\left(\sum_{\sigma;|\sigma|+\theta>|k|}({\sf\Pi}_x^{\sf g}\sigma){\sf g}_x(\tau/\sigma)\right)(x).
\end{align*}
\end{Dem}

\medskip

\begin{defn}
\label{DefnAdmissibleModel}
A \textbf{\textsf{model}} $\sf (g,\Pi)$ on $\mathscr{T}$ is said to be \textbf{\textsf{admissible}} if the identities ${\sf g}_x(\bsX_e^k)={\sf\Pi}\bsX_e^k(x)=x_e^k$, the commutation rule \eqref{commutation I and Pi}, and \eqref{EqDefnAdmissibleGBis} are satisfied.
\end{defn}

\medskip

\begin{Remark*}
Note that our notion of admissible model is more general than the corresponding notion introduced by Bruned, Hairer and Zambotti in \cite{BHZ16}; Definition 6.9 in that work. Their admissible $\sf\Pi$-maps, together with the positive twisted antipode from Proposition 6.2 in \cite{BHZ16}, are used in definition 6.8 therein to build a $\sf g$-map and models $\sf (g,\Pi)$ that are admissible model on $\mathscr{T}$ in our sense, with all $\phi_e\equiv 1$. This is a direct consequence of Lemma 6.10 in \cite{BHZ16} and the following equalities.
\begin{align}
\label{admissible Pi(It)}{\sf\Pi}_x^{\sf g}\mcI\tau
&=\bfI({\sf\Pi}_x^{\sf g}\tau)-
\sum_{e\in E, |k|<|\tau|+\theta} \frac{((\cdot)_e-x_e)^k}{k!}\,\bfI_k^e({\sf\Pi}^{\sf g}_x\tau)(x),\\
\label{admissible g(It)}{\sf g}_x^{-1}(\mcI_k^{e+}\tau)
&=-\sum_{\ell;|k+\ell|<|\tau|+\theta}\frac{(-x_e)^\ell}{\ell!}\bfI_{k+\ell}^e({\sf\Pi}_x^{\sf g}\tau)(x),\\
\label{admissible g_yx(It)}{\sf g}_{yx}(\mcI_k^{e+}\tau)
&=\sum_{\sigma\le\tau,|k|<|\sigma|+\theta}{\sf g}_{yx}(\tau/\sigma)\bfI_k^e({\sf\Pi}_y^{\sf g}\sigma)(y)
-\sum_{|k+\ell|<|\tau|+\theta}\frac{(y_e-x_e)^\ell}{\ell!}\bfI_{k+\ell}^e({\sf\Pi}_x^{\sf g}\tau)(x).
\end{align}
Let us show the above equalities. First we assume that \eqref{admissible Pi(It)} holds for any $\sigma\in\mcB$ with $\sigma<\tau$. Then by \eqref{commutation I and Pi} and \eqref{EqDefnAdmissibleGBis},
\begin{align*}
{\sf\Pi}_x^{\sf g}\mcI\tau
&={\sf\Pi}\mcI\tau-\sum_{\sigma<\tau}{\sf g}_x(\tau/\sigma){\sf\Pi}_x^{\sf g}\mcI\sigma-\sum_{e,k}\frac{((\cdot)_e-x_e)^k}{k!}{\sf g}_x(\mcI_k^{e+}\tau)\\
&=\bfI({\sf\Pi}\tau)-\sum_{\sigma<\tau}{\sf g}_x(\tau/\sigma)\left(\bfI({\sf\Pi}_x^{\sf g}\sigma) -
\sum_{e\in E, |k|<|\sigma|+\theta} \frac{((\cdot)_e-x_e)^k}{k!}\,\bfI_k^e({\sf\Pi}^{\sf g}_x\sigma)(x)\right)   \\
&\quad-\sum_{e,|k|<|\tau|+\theta}\frac{((\cdot)_e-x_e)^k}{k!}
\left(\bfI_k^e\big({\sf\Pi}^{\sf g}_x\tau\big)(x) + \sum_{\sigma<\tau, |k|<|\sigma|+\theta}\, {\sf g}_x(\tau/\sigma)\, \bfI_k^e\big({\sf\Pi}^{\sf g}_x\sigma\big)(x)\right)   \\
&=\bfI\big({\sf\Pi}\tau-\sum_{\sigma<\tau}{\sf g}_x(\tau/\sigma){\sf\Pi}_x^{\sf g}\sigma\big)
-\sum_{e,|k|<\tau+\theta}\frac{((\cdot)_e-x_e)^k}{k!}\bfI_k^e\big({\sf\Pi}^{\sf g}_x\tau\big)(x),
\end{align*}
where by using ${\sf\Pi}\tau-\sum_{\sigma<\tau}{\sf g}_x(\tau/\sigma){\sf\Pi}_x^{\sf g}={\sf\Pi}_x^{\sf g}\tau$ we have \eqref{admissible Pi(It)}. On the other hand, by \eqref{EqIntertwiningDeltaDelta+},
\begin{align*}
0={\sf g}_{xx}(\mcI_{k+\ell}^{e+}\tau)
&=\sum_{\sigma\le\tau}{\sf g}_x(\mcI_{k+\ell}^{e+}\sigma){\sf g}_x^{-1}(\tau/\sigma)
+\sum_{m;|k+\ell+m|<|\tau|+\theta}\frac{x_e^m}{m!}{\sf g}_x^{-1}(\mcI_{k+\ell+m}^{e+}\tau)
\end{align*}
From this identity and \eqref{EqDefnAdmissibleGBis},
\begin{align*}
{\sf g}_x^{-1}(\mcI_k^{e+}\tau)
&=\sum_{\ell,m}\frac{(-x_e)^\ell}{\ell!}\frac{x_e^m}{m!}{\sf g}_x^{-1}(\mcI_{k+\ell+m}^{e+}\tau) = -\sum_{\sigma\le\tau,\ell}\frac{(-x_e)^\ell}{\ell!}{\sf g}_x(\mcI_{k+\ell}^{e+}\sigma){\sf g}_x^{-1}(\tau/\sigma)\\
&=-\sum_{\sigma\le\tau,\ell}\frac{(-x_e)^\ell}{\ell!}{\sf g}_x^{-1}(\tau/\sigma)
\sum_{\eta\le\sigma, |k+\ell|<|\eta|+\theta}\, {\sf g}_x(\sigma/\eta)\, \bfI_{k+\ell}^e\big({\sf\Pi}^{\sf g}_x\eta\big)(x)\\
&=-\sum_{\eta\le\tau,|k+\ell|<|\eta|+\theta}\frac{(-x_e)^\ell}{\ell!}{\sf g}_{xx}(\tau/\eta)\bfI_{k+\ell}^e({\sf\Pi}_x^{\sf g}\eta)(x),
\end{align*}
where ${\sf g}_{xx}(\tau/\eta)=\mathbf{1}_{\eta=\tau}$, which yields \eqref{admissible g(It)}.
Moreover, we have \eqref{admissible g_yx(It)} as follows.
\begin{align*}
{\sf g}_{yx}(\mcI_k^{e+}\tau)
&=\sum_{\sigma\le\tau}{\sf g}_y(\mcI_k^{e+}\sigma){\sf g}_x^{-1}(\tau/\sigma)+\sum_\ell\frac{y_e^\ell}{\ell!}{\sf g}_x^{-1}(\mcI_{k+\ell}^{e+}\tau)   \\
&=\sum_{\eta\le\sigma\le\tau,|k|<|\eta|+\theta}{\sf g}_y(\sigma/\eta){\sf g}_x^{-1}(\tau/\sigma)\bfI_k^e({\sf\Pi}_y^{\sf g}\eta)(y)   \\
&\quad- \sum_{\ell,m;|k+\ell+m|<|\tau|+\theta}\frac{y_e^\ell}{\ell!}\frac{(-x_e)^m}{m!}\bfI_{k+\ell+m}^e({\sf\Pi}_x^{\sf g}\tau)(x)   \\
&=\sum_{\eta\le\tau,|k|<|\eta|+\theta}{\sf g}_{yx}(\tau/\eta)\bfI_k^e({\sf\Pi}_y^{\sf g}\eta)(y) - \sum_{|k+\ell|<|\tau|+\theta}\frac{(y_e-x_e)^\ell}{\ell!}\bfI_{k+\ell}^e({\sf\Pi}_x^{\sf g}\tau)(x).
\end{align*}
\end{Remark*}

\medskip

\subsection{Parametrization of the set of admissible models}

We prove Theorem \ref{ThmMain2} in this section, giving a parametrization of the set of admissible models by a product of H\"older spaces of non-positive regularity exponent. A similar parametrization of the space of branched rough paths was achieved in the recent work \cite{TZ18} of Tapia and Zambotti, with very different tools. The next result applies in particular in the former setting, when formulated in terms of regularity structures. We need the following structural assumptions on $T^+$ and $T$.

\medskip

\noindent \textsf{\textbf{Assumption (C)}}   \vspace{0.1cm}
{\it
\begin{itemize}
\setlength{\parskip}{0.1cm}
\item The basis $\mcB^+$ of $T^+$ is a commutative monoid with unit $\bf 1$, freely generated by the set
$$
\{X_e^i\}_{e\in E, i=1,\dots,d} \cup \big\{\mcI_k^{e+}\tau\big\}_{\tau\in\mcB, e\in E, k\in\bbN^d,|\tau|+\theta-|k|>0},
$$
\item For any $\tau,\sigma\in\mcB$, the element $\tau/\sigma\in T^+$ is contained in the subalgebra generated by the set
$$
\{X_e^i\}_{e\in E, i=1,\dots,d} \cup \big\{\mcI_k^{e+}\eta\big\}_{\eta\in\mcB, e\in E, k\in\bbN^d; |\eta|<|\tau|, |\eta|+\theta-|k|>0},
$$
\end{itemize}   }

\medskip

\begin{thm}
\label{ThmModelsFromBrackets}
Let  a regularity structure $\mathscr{T}$ equipped with an abstract integration map satisfy assumptions \textbf{\textsf{(A)}}, \textbf{\textsf{(B)}} and \textbf{\textsf{(C)}}. Given any family of distributions $\big(\Brac{\tau}\in\mcC^{|\tau|}(\bbR^d)\big)_{\tau\in\mcB; |\tau|\leq 0}$, there exists a unique admissible model $\sf (g,\Pi)$ on $\mathscr{T}$ such that one has 
\begin{equation}
\label{EqStructurePiBracket}
{\sf \Pi}\tau := \sum_{\sigma<\tau}{\sf P}_{{\sf g}(\tau/\sigma)}\Brac{\sigma} + \Brac{\tau},
\end{equation}
for all $\tau\in\mcB$ with $|\tau|\leq 0$.
\end{thm}

\medskip

\begin{Dem}
For $\alpha\in A$, define $T_{(\alpha)}^+$ as the subalgebra of $T^+$ generated by
$$
\{X_e^i\}_{e\in E, i=1..d} \cup \big\{\mcI_k^{e+}\tau\big\}_{\tau\in\mcB, e\in E, k\in\bbN^d,|\tau|<\alpha}.
$$
By Assumption \textbf{\textsf{(C)}}, $T_{(\alpha)}^+$ is closed under $\Delta^+$. Start by noting that $\mathscr{T}_{<\alpha} := (T_{(\alpha)}^+,T_{<\alpha})$ is a regularity structure for any $\alpha\in A$. Define inductively on $\alpha\in A$ the maps 
$$
{\sf\Pi}_{<\alpha} : T_{<\alpha} \mapsto \mcC^{\beta_0}(\bbR^d),
$$
and 
$$
{\sf g}^{(\alpha)} : T_{(\alpha)}^+ \mapsto C_b(\bbR^d),
$$
with ${\sf g}_x^{(\alpha)}(\bsX_e^k) = x_e^k$, initializing the induction. Write ${\sf M}_{<\alpha}$ for the model $({\sf g}^{(\alpha)}, {\sf \Pi}_{<\alpha})$ on $\mathscr{T}_{<\alpha}$. Set $\alpha^+ := \min\{\beta>\alpha; \beta\in A\}$. Given a basis vector $\tau\in \mcB_{\alpha}$, the function $\bsh_\tau := \sum_{\sigma<\tau} {\sf g}^{(\alpha)}(\tau/\sigma)\,\sigma$ is an element of $\mcD^{\alpha}(\mathscr{T}_{<\alpha},{\sf g}^{(\alpha)})$. Define ${\sf \Pi}_{<\alpha^+}\tau$ as equal to either
$$
\sum_{\sigma<\tau} {\sf P}_{{\sf g}^{(\alpha)}(\tau/\sigma)}\Brac{\sigma} + \Brac{\tau},
$$
if $|\tau|\leq 0$, or
$$
\textsf{\textbf{R}}^{{\sf M}_{<\alpha}}(\bsh_\tau),
$$
if $|\tau|>0$, where $\textsf{\textbf{R}}^{{\sf M}_{<\alpha}}$ stands for the reconstruction operator on $\mcD^{\alpha}\big(\mathscr{T}_{<\alpha},{\sf g}^{(\alpha)}\big)$ associated with the model ${\sf M}_{<\alpha}$. We have in both cases $\big|\big\langle ({\sf \Pi}_{<\alpha^+})_x^{{\sf g}^{(\alpha)}}\tau ,  \varphi_x^\lambda\big\rangle\big| \lesssim \lambda^{\alpha}$, from Corollary \ref{CorUniqueExtension}. Define then an extension ${\sf g}^{(\alpha^+)}$ of ${\sf g}^{(\alpha)}$ onto $T_{(\alpha^+)}^+$ by requiring that it is multiplicative and by defining ${\sf g}^{(\alpha^+)}(\mcI_k^{e+}\tau)$ from identity \eqref{EqDefnAdmissibleGBis}, with ${\sf\Pi}_{<\alpha^+}$ in the role of $\sf\Pi$. Boundedness of ${\sf g}^{(\alpha^+)}$ is checked by induction. Given Assumption \textbf{\textsf{(C)}} on the regularity structure $\mathscr{T}$, closing the induction step amounts to proving that 
$$
\big|{\sf g}^{(\alpha^+)}_{yx}(\mcI^{e+}_k\tau)\big| \lesssim |y-x|^{|\tau| + \theta - |k|},
$$ 
for every $k\in\bbN^d$ and $e\in E$. Set
$$
\big(\Upsilon_\alpha^e\bsh_\tau\big)(x) := \mcI_0^{e+}\bsh_\tau(x) + \sum_{k\in\bbN^d} {\sf g}_x(\mcI^{e+}_k\tau)\,\frac{\bsX_e^k}{k!}.
$$
Proposition \ref{PropRealizationPiPi+} is used to prove the following fact, proved below.

\medskip

\begin{lem}
\label{LemConstructionAdmissibleModels}
One has $\Upsilon_\alpha^e\bsh_\tau\in\mcD^{\alpha+\theta}\big(\mathscr{T}_{<\alpha}^+, {\sf g}^{(\alpha)}\big)$.
\end{lem}

\medskip

But 
$$
\big(\Upsilon_\alpha^e\bsh_\tau\big)(y) - \widehat{{\sf g}^{(\alpha)}_{yx}}\big(\Upsilon_\alpha^e\bsh_\tau\big)(x)
$$ 
has $\bsX_e^k$ component equal to 
$$
\left| {\sf g}_y(\mcI_\ell^{e+}\tau)  -  \sum_{\eta<\tau} {\sf g}_x(\tau/\eta){\sf g}_{yx}(\mcI_\ell^{e+}\eta)  -  \sum_m {\sf g}_x(\mcI_{\ell+m}^{e+}\tau)\frac{(y_e - x_e)^m}{m!} \right|  \lesssim  |y-x|^{|\tau|+\theta-|\ell|},
$$
from lemma \ref{LemConstructionAdmissibleModels}; we recognize $\big|{\sf g}^{(\alpha^+)}_{yx}(\mcI^{e+}_k\tau)\big|$ in the left hand side, which closes the induction.
\end{Dem}

\smallskip

\noindent\textbf{\textsf{Proof of Lemma \ref{LemConstructionAdmissibleModels} -- }} We follow closely the proof of Schauder estimates for modelled distributions -- Theorem 5.12 of \cite{Hai14}. Note first that by \eqref{EqDefnAdmissibleGBis} one can decompose $\Upsilon_\alpha^e\bsh_\tau$ under the form
$$
\big(\Upsilon_\alpha^e\bsh_\tau\big)(x) = \mcI_0^{e+}\bsh_\tau(x) + \mcJ^e(x)(\bsh_\tau(x)) + (\mcN^e\bsh_\tau)(x),
$$
with 
$$
\mcJ^e(x)\bsh_\tau(x) := \sum_{\sigma<\tau} {\sf g}_x(\tau/\sigma) \sum_{|k|<|\sigma|+\theta} \frac{\bsX_e^k}{k!}\, \bfI_k^e({\sf\Pi}_x^{\sf g}\sigma)(x)
$$
and 
\begin{equation*}
\big(\mcN^e\bsh_\tau\big)(x) 
:= \sum_{|k|<|\tau|+\theta} \frac{\bsX_e^k}{k!} \, \bfI_k^e\big(\textsf{\textbf{R}}(\bsh_\tau)-{\sf\Pi}_x^{\sf g}\bsh_\tau(x)\big)(x)\\
= \sum_{|k|<|\tau|+\theta} \frac{\bsX_e^k}{k!} \, \bfI_k^e\big({\sf\Pi}_x^{\sf g}\tau\big)(x),
\end{equation*}
with $|\tau|=\alpha$, where $\mcJ^e(x)$ is an operator on $T$ rather than on $\mcD^\gamma(\mathscr{T},{\sf g})$, defined by
$$
\mcJ^e(x)\sigma=\sum_{|k|<|\sigma|+\theta}\frac{\bsX_e^k}{k!} \bfI_k^e({\sf\Pi}^{\sf g}_x\sigma)(x).
$$
Remark then, as in Lemma 5.16 of \cite{Hai14}, that we have for any $x,y\in\bbR^d$ 
\begin{equation}
\label{EqCommutationRelation}
\widehat{{\sf g}_{yx}}^+\big(\mcI^{e+}+\mcJ^e(x)\big) = \big(\mcI^{e+}+\mcJ^e(y)\big)\widehat{{\sf g}_{yx}}.
\end{equation}
We give a direct proof. By definition, we have
\begin{align*}
\widehat{{\sf g}_{yx}}^+\big(\mcI^{e+}+\mcJ^e(x)\big) \tau &= \sum_\sigma{\sf g}_{yx}(\tau/\sigma)\mcI^{e+}\sigma + \sum_{k}\frac{\bsX_e^k}{k!}{\sf g}_{yx}(\mcI_{k}^{e+}\tau)   \\
&\quad+\sum_{k,\ell;|k+\ell|<|\tau|+\theta}\frac{\bsX_e^k}{k!}\frac{(y_e-x_e)^\ell}{\ell!}\bfI_{k+\ell}^e({\sf\Pi}_x^{\sf g}\tau)(x)
\end{align*}
and
\begin{align*}
(\mcI^{e+}+\mcJ^e(y))\widehat{{\sf g}_{yx}}\tau
=\sum_\sigma{\sf g}_{yx}(\tau/\sigma)\mcI^{e+}\sigma
+\sum_{k,\sigma;|k|<|\sigma|+\theta}\frac{\bsX_e^k}{k!}
{\sf g}_{yx}(\tau/\sigma)\bfI_k^e({\sf\Pi}_y^{\sf g}\sigma)(y).
\end{align*}
They are equal because of \eqref{admissible g_yx(It)}. We use the interwining relation \eqref{EqCommutationRelation} to write
\begin{align*}
\big(\Upsilon_\alpha^e&\bsh_\tau\big)(y) - \widehat{{\sf g}_{yx}}^+\big(\Upsilon_\alpha^e\bsh_\tau\big)(x) 
= \big(\Upsilon_\alpha^e\bsh_\tau\big)(y) - \widehat{{\sf g}_{yx}}^+ \big(\mcI^{e+}+\mcJ^e(x)\big)\bsh_\tau(x) - \widehat{{\sf g}_{yx}}^+(\mcN^e\bsh_\tau)(x)   \\
&= \big(\Upsilon_\alpha^e\bsh_\tau\big)(y) - \big(\mcI^{e+}+\mcJ^e(y)\big)\widehat{{\sf g}_{yx}}\bsh_\tau(x) - \widehat{{\sf g}_{yx}}^+(\mcN^e\bsh_\tau)(x)   \\
&= \mcI^{e+}\Big(\bsh_\tau(y) - \widehat{{\sf g}_{yx}}\bsh_\tau(x)\Big) + \mcJ^e(y)\Big(\bsh_\tau(y) - \widehat{{\sf g}_{yx}}\bsh_\tau(x)\Big) + \Big((\mcN^e\bsh_\tau)(y) - \widehat{{\sf g}_{yx}}^+(\mcN^e\bsh_\tau)(x)\Big).
\end{align*}
For the $\mcI^{e+}$ term, one has the elementary estimate
$$
\big\|\mcI^{e+}\big(\bsh_\tau(y) - \widehat{{\sf g}_{yx}}\bsh_\tau(x)\big)\big\|_\beta
\lesssim \big\| \bsh_\tau(y) - \widehat{{\sf g}_{yx}}\bsh_\tau(x) \big\|_{\beta-\theta}
\leq \|\bsh_\tau\|_{\mcD^\alpha} \, |y-x|^{\alpha+\theta-\beta}.
$$
The $\mcJ^e$ and $\mcN^e$ terms take values in the polynomial part of $T^+$. Write $\tau^{\bsX_e^k}$ for the $\bsX_e^k$-component of $\tau\in T^+$. Decompose $K(y,z)=\sum_{n=0}^\infty K_n(y,z)$, and let $\admi^e=:\sum_n\admi_n^e$ and $\remind^e = :\sum_n\remind_n^e$, be the corresponding operators. We have
\begin{align*}
&\Big(\mcJ_n^e(y)\big(\bsh_\tau(y) - \widehat{{\sf g}_{yx}}\bsh_\tau(x)\big) + \big(\mcN_n^e\bsh_\tau)(y) - \widehat{{\sf g}_{yx}}^+(\mcN_n^e\bsh_\tau)(x)\Big)^{\bsX_e^k}\\
&=\Big(\mcJ_n^e(y)\big(\bsh_\tau(y) - \widehat{{\sf g}_{yx}}\bsh_\tau(x)\big) \Big)^{\bsX_e^k}+ \Big(\big(\mcN_n^e\bsh_\tau)(y) - \widehat{{\sf g}_{yx}}^+(\mcN_n^e\bsh_\tau)(x)\Big)^{\bsX_e^k}
=: (\Asterisk)^1_n + (\Asterisk)^2_n
\end{align*}
and
\begin{align*}
&\Big(\mcJ_n^e(y)\big(\bsh_\tau(y) - \widehat{{\sf g}_{yx}}\bsh_\tau(x)\big) + \big(\mcN_n^e\bsh_\tau)(y) - \widehat{{\sf g}_{yx}}^+(\mcN_n^e\bsh_\tau)(x)\Big)^{\bsX_e^k}\\
&= \sum_{\beta\in A, |k|<\beta+\theta}
\int_{\bbR^d}\frac1{k!}\partial_y^k\big(\phi_e(y)K_n(y,z)\big) \, {\sf\Pi}^{\sf g}_y \big(\bsh_\tau(y) - \widehat{{\sf g}_{yx}}\bsh_\tau(x)\big)_\beta(z)\,dz   \\
&\quad+ \int_{\bbR^d}\frac1{k!}\partial^k(\phi_e K_n)_{y,x}^{\alpha+\theta-|k|}(z)\, ({\sf\Pi}^{\sf g}_x\tau)(z)\,dz     \\
&\quad + \int_{\bbR^d}\frac1{k!}\partial_y^k(\phi_e(y)K_n(y,z))\,{\sf\Pi}^{\sf g}_y\big(\widehat{{\sf g}_{yx}}\bsh_\tau(x)-\bsh_\tau(y)\big)(z) \, dz  \\
&=\int_{\bbR^d}\frac1{k!}\partial^k(\phi_e K_n)_{y,x}^{\alpha+\theta-|k|}(z)\,({\sf\Pi}^{\sf g}_x\tau)(z)\,dz   \\
&\quad+ \sum_{\beta\in A,|k|\ge\beta+\theta} \int_{\bbR^d}\frac1{k!}\partial_y^k\big(\phi_e(y)K_n(y,z)\big) \,{\sf\Pi}^{\sf g}_y \big(\widehat{{\sf g}_{yx}}\bsh_\tau(x)-\bsh_\tau(y)\big)_\beta(z) \, dz   \\
&=:(\bigstar)_n^1+(\bigstar)_n^2,
\end{align*}
where
\begin{align*}
\partial^k(\phi_e K_n)_{y,x}^{\alpha+\theta-|k|}(z)
:=\partial_y^k\big(\phi_e(y) K_n(y,z)\big) - \sum_{|\ell|<\alpha+\theta-|k|}\frac{(y_e-x_e)^\ell}{\ell!}\,\partial_x^{k+\ell}\big(\phi_e(x)K_n(x,z)\big).
\end{align*}
Write $|y-x|\simeq 2^{-N}$. We use the $(\Asterisk)$-decomposition, with $\mcJ^n$ and $\mcN^n$ separated, to estimate the sum over $n>N$, and the $(\bigstar)$-decomposition to estimate the sum over $n\leq N$. Each decomposition is well-adapted to get $N$-independent upper bounds. 

\smallskip

$\bullet$ For $n>N$, we have from the bound \eqref{Bound zeta_x(K_n)} and its derivatives the estimate
\begin{align*}
\sum_{n>N} \big|(\Asterisk)^1_n\big| \lesssim  \sum_{n>N} \sum_{\beta\in A,|k|<\beta+\theta} |y-x|^{\alpha-\beta}\,2^{(|k|-\beta-\theta)n}  \lesssim  |y-x|^{\alpha+\theta-|k|}.
\end{align*}
By the definition of $\mcN$, we get
\begin{align*}
\sum_{n>N} \big|(\Asterisk)^2_n\big|
&\leq \sum_{n>N} \left| \int_{\bbR^d}\frac1{k!}\partial_y^k\big(\phi_e(y)K_n(y,z)\big) ({\sf\Pi}^{\sf g}_y\tau)(z)dz\right|   \\
&\quad+ \sum_{n>N} \left| \sum_{|\ell|+|k|<\alpha+\theta}\frac{(y_e-x_e)^\ell}{\ell!}\int_{\bbR^d}\frac1{k!}\partial_x^{k+\ell}\big(\phi_e(x)K_n(x,z)\big)({\sf\Pi}^{\sf g}_x\tau)(z)dz \right|   \\
&\lesssim  \sum_{n>N} \left(2^{(|k|-\alpha-\theta)n}+\sum_{|\ell|<\alpha+\theta-|k|}|y-x|^\ell 2^{(|\ell|-\alpha-\theta)n}\right) \lesssim |y-x|^{\alpha+\theta-|k|}
\end{align*}

\smallskip

$\bullet$ To deal with the sum over $n\leq N$, we use the $(\bigstar)$-decomposition. For $(\bigstar)_1^n$, note that since $|y-x|\simeq 2^{-N}\le 2^{-n}$, the function $\partial^k(\phi_e K_n)_{y,x}^{\alpha+\theta-|k|}$ is supported on a ball $B(x,C2^{-n})$, for some positive constant $C$. From Taylor formula with bounded polynomials proved in Appendix \ref{SectionBoundedPolynomials}, we have 
\begin{align}\label{Taylor phi_eK_n}
\Big|\partial_z^m\partial^k(\phi_e K_n)_{y,x}^{\alpha+\theta-|k|}(z)\Big|
&\lesssim B_{\alpha+\theta} (\partial_z^m K_n (\cdot,z)) \, |y-x|^{\alpha+\theta-|k|},
\end{align}
with $B_r\big(\partial_z^m K_n(\cdot,z)\big)\lesssim 2^{(d+|m|+r-\theta-\epsilon) n}$, from either \eqref{EqConditionKn} or the interpolation theorem 2.80 in \cite{BCD11}. Hence
\begin{align*}
\Big|\partial_z^m\partial^k(\phi_e K_n)_{y,x}^{\alpha+\theta-|k|}(z)\Big| &\lesssim2^{(d+|m|+\alpha-\epsilon)n}|y-x|^{\alpha+\theta-|k|}.
\end{align*}
It follows then from the proof of Lemma \ref{LemmeDefnIDistribution} that we have
\begin{align*}
|(\bigstar)_n^1|&= \Big|({\sf\Pi}_x^{\sf g}\tau)\big(\partial^k(\phi_e K_n)_{y,x}^{\alpha+\theta-|k|}\big)\Big| \lesssim 2^{-\epsilon n}|y-x|^{\alpha+\theta-|k|},
\end{align*}
so the sum over $n\leq N$ is independent of $N$, of order $|y-x|^{\alpha+\theta-|k|}$. As for the $(\bigstar)_2^n$-terms, they involve some indices $\zeta$ with $|k|\ge\zeta+\theta$, so the same elementary bounds as above give
$$
|(\bigstar)_n^2| \lesssim \sum_{\zeta\in A,|k|\ge\zeta+\theta} |y-x|^{\alpha-\zeta} \, 2^{(|k|-\zeta-\theta-\epsilon)n} \lesssim 2^{-\epsilon n} \, |y-x|^{\alpha+\theta-|k|},
$$
since $2^n\leq |y-x|^{-1}$. The sum over $n\leq N$ of the $(\bigstar)_2^n$ is thus independent of $N$, of order $|y-x|^{\alpha+\theta-|k|}$.   \hfill$\rhd$

\medskip

\begin{Remark*}
If $(T^+,T)$ satisfies $|\tau|+\theta\notin\bbN$ for any $\tau\in\mcB$, then we can choose $\varepsilon=0$ for the estimate \eqref{EqConditionKn} on the kernel $K_n$. We need to modify the argument for the sum over $n\le N$. For $(\bigstar)_n^1$, since $\partial^k(\phi_e K_n)_{y,x}^{\alpha+\theta-|k|}=\partial^k(\phi_e K_n)_{y,x}^{\alpha+\theta+\delta-|k|}$ in \eqref{Taylor phi_eK_n} for small $\delta>0$ such that $(\alpha+\theta,\alpha+\theta+\delta)\cap\bbN=\emptyset$, we have
$$
|(\bigstar)_n^1|\lesssim2^{\delta n}|y-x|^{\alpha+\theta+\delta-|k|},
$$
so the sum over $n\le N$ is of order $|y-x|^{\alpha+\theta-|k|}$. For $(\bigstar)_n^2$, since they involve indices $\zeta$ with $|k|>\zeta+\theta$, we have
$$
\sum_{n\le N}|(\bigstar)_n^2|\lesssim
\sum_{\zeta\in A,|k|>\zeta+\theta}|y-x|^{\alpha-\zeta}2^{(|k|-\zeta-\theta)N}
\lesssim|y-x|^{\alpha+\theta-|k|}.
$$
\end{Remark*}

\bigskip

\appendix

\section{Paraproducts}
\label{AppendixParapdcts}

We summarize in this section some basic concepts and results of the Littlewood-Paley theory. Let $\{\rho_i\}_{i=-1}^\infty$ be a dyadic partition of unity of $\Euc^d$, i.e., $\rho_i:\Euc^d\to[0,1]$ is a compactly supported smooth radial function with the following properties.
\begin{itemize}
\item $\supp{\rho_{-1}}\subset\big\{x\in\Euc^d;|x|<\frac43\big\}$ and $\supp{\rho_0}\subset\big\{x\in\Euc^d;\frac34<|x|<\frac83\big\}$.
\item $\rho_i(x)=\rho_0(2^{-i}x)$ for any $x\in\Euc^d$ and $i\ge0$.
\item $\sum_{i=-1}^\infty\rho_i(x)=1$ for any $x\in\Euc^d$.
\end{itemize}
We define the \emph{Littlewood-Paley blocks} $\{\Delta_i\}_{i=-1}^\infty$ by
$$
\Delta_if:=\Fouri^{-1}(\rho_i\Fouri f),\quad f\in\Schw'(\Euc^d),
$$
where $\Fouri$ is a Fourier transform on $\Euc^d$ defined by
$$
\Fouri\varphi(\xi):=\int_{\Euc^d}\varphi(x)e^{-2\pi\imunit\innpro{x}{\xi}}dx,\quad\varphi\in\Schw(\Euc^d).
$$
Now we define the \emph{H\"older-Besov spaces}. For any $\alpha\in\Euc$ and $f\in\Schw'(\Euc^d)$, we define
$$
\|f\|_{\Hoel^\alpha}:=\sup_{i\ge-1}2^{\alpha i}\|\Delta_if\|_{L^\infty(\Euc^d)}.
$$
We denote by $\mcC^\alpha(\bbR^d)$ the space of all $f\in\mcS'(\bbR^d)$ with $\|f\|_{\mcC^\alpha}<\infty$. This definition does not ensure the separability of $\Hoel^\alpha(\bbR^d)$, so it may be better to consider the space $\Hoel^\beta_0(\bbR^d)$, the completion of $\Schw(\Euc^d)$ under the norm $\|\cdot\|_{\Hoel^\alpha}$. However, it does not matter because $\Hoel^\alpha(\bbR^d)$ is embedded into the space $\Hoel^\beta_0(\bbR^d)$ for any $\beta<\alpha$ -- see e.g. \cite[Proposition~2.74]{BCD11}. For $\alpha\in(0,\infty)\setminus\Nat$, the norm $\|f\|_{\Hoel^\alpha}$ is equivalent to the H\"older norm
$$
\|f\|_{C^\alpha}:=\sum_{k\in\Nat^d,|k|<\alpha}\|\partial^kf\|_{L^\infty}+\sum_{k\in\Nat^d,|k|=\lceil\alpha\rceil}\|\partial^kf\|_{\text{($\alpha-\lceil\alpha\rceil$)-H\"older}},
$$
where $\|\partial^kf\|_{\text{($\alpha-\lceil\alpha\rceil$)-H\"older}}$ is the infimum of constants $C$ such that the property 
$$
|\partial^kf(y)-\partial^kf(x)|\le C|y-x|^{\alpha-\lceil\alpha\rceil},
$$ 
holds for any $x,y\in\Euc^d$.
For $\alpha\in\Nat$, the space $\Hoel^\alpha(\bbR^d)$ is strictly larger than the space $C_b^\alpha(\bbR^d)$ with the norm
$$
\|f\|_{C_b^\alpha}:=\sum_{k\in\Nat^d,|k|\le\alpha}\|\partial^kf\|_{L^\infty};
$$
see e.g. \cite[page~99]{BCD11}. Bony's \emph{paraproduct} ${\sf P}$ and \emph{resonant} operator $\circleddash$ are defined by
$$
{\sf P}_fg := \sum_{\substack{i,j\ge-1\\ i\le j-2}}\Delta_if\Delta_jg,\quad
\circleddash(f,g) := \sum_{\substack{i,j\ge-1\\|i-j|\le1}}\Delta_if\Delta_jg,
$$
for any $f,g\in\Schw'(\Euc^d)$, as long as they converge. We then have formally 
$$
fg={\sf P}_fg+{\sf P}_gf+\circleddash(f,g).
$$
The basic continuity results for these operators are summarized as follows.

\smallskip

\begin{prop}[{\cite[Theorems~2.82 and 2.85]{BCD11}}]
Let $\alpha,\beta\in\Euc$.   \vspace{0.1cm}

\begin{enumerate}
	\item $\|{\sf P}_fg\|_{\Hoel^\beta}\lesssim\|f\|_{L^\infty}\|g\|_{\Hoel^\beta}$.   \vspace{0.1cm}
	
	\item If $\alpha<0$, then $\|{\sf P}_fg\|_{\Hoel^{\alpha+\beta}}\lesssim\|f\|_{\Hoel^\alpha}\|g\|_{\Hoel^\beta}$.   \vspace{0.1cm}

	\item If $\alpha+\beta>0$, then $\big\|\circleddash(f,g)\big\|_{\Hoel^{\alpha+\beta}}\lesssim\|f\|_{\Hoel^\alpha}\|g\|_{\Hoel^\beta}$.
\end{enumerate}
\end{prop}

\bigskip

\section{Bounded polynomials}
\label{SectionBoundedPolynomials}

This appendix is a follow-up of example 1 in Section \ref{SectionBasicsRS} describing bounded polynomials and their associated regularity structure. We give the proofs of Proposition \ref{app: prop PTE} and Proposition \ref{PropPolynomialLift}. Set $\Lambda := (\frac14\Int)^d$, and, for any $\lambda = (\lambda_i)_{i=1..d}\in \Lambda$, define $U_\lambda := \prod_{i=1}^d\big(\lambda_i-\frac3{16},\lambda_i+\frac3{16}\big)$. This family of bounded open subsets of $\bbR^d$ cover $\bbR^d$, and are uniformly locally finite covering, i.e.,
$$
\sup_{x\in\Euc^d}\#\big\{\lambda\in\Lambda;x\in U_\lambda\big\}<\infty.
$$
For $x\in\bbR^d$ and $A\subset\bbR^d$, set $d(x,A) := \inf\big\{|x-y|;y\in A\big\}$.

\medskip

\begin{lem*}
One can construct tuples $\big(\varphi_\lambda, (\psi_\lambda^i)_{i=1..d}\big)_{\lambda\in\Lambda}$ of smooth real-valued functions on $\bbR^d$, with compact support, with the following properties.   \vspace{0.15cm}

\begin{enumerate}
	\item[($\textsf{a}_1$)] One has $\varphi_\lambda(x)\ge0$, for any $x\in\Euc^d$, and $\varphi_\lambda(x)=0$ for any $x\in U_\lambda^c$, for any $\lambda\in\Lambda$.

	\item[($\textsf{a}_2$)] One has $\sum_{\lambda\in\Lambda} \varphi_\lambda(x)=1$, for any $x\in\Euc^d$.

	\item[($\textsf{a}_3$)] For any $N>0$, there is a constant $C_N$ independent of $\alpha$ such that 
	$$
	\big|(\partial^\ell\varphi_\lambda)(x)\big|\le C_N\,d(x,U_\lambda^c)^N, 
	$$
	for any $\lambda\in\Lambda, x\in\Euc^d$, and $|\ell|\le N$.   \vspace{0.15cm}
\end{enumerate}

\begin{enumerate}
	\item[($\textsf{b}$)] The functions $\psi_\lambda^i$ are uniformly bounded and $\psi_\lambda^i(y)-\psi_\lambda^i(x)=y_i-x_i$, for any $x,y\in U_\lambda$.
\end{enumerate}
\end{lem*}

\medskip

\begin{Dem}
We let the reader construct a partition of unity $\{\varphi_\lambda\}$ satisfying assumptions $\textsf{a}_1$ to $\textsf{a}_3$. The third property is ensured if we impose
$$
\varphi_\lambda(x)\simeq \exp\left(-\frac1{\big|x_i-(\lambda_i\pm\frac3{16})\big|}\right)
$$
when $x\in U_\lambda$ is near $\partial U_\lambda\cap\left\{x_i=\lambda_i\pm\frac3{16}\right\}$.
For each $\lambda\in\Lambda$, we choose a smooth function $\psi_\lambda^i$ such that
$$
\psi_\lambda^i(x)=
\begin{cases}
x_i-\lambda_i&x\in U_\lambda,\\
0&x\notin V_\lambda:=\prod_{i=1}^d\big(\lambda_i-\frac14,\lambda_i+\frac14\big),
\end{cases}
$$
and $|\psi_\lambda^i(x)|\le1$, for any $x\in\Euc^d$.
\end{Dem}

\medskip

Recall we define $B_r(f):=\|f\|_{C_b^r}$ if $r\in\Nat$, and $B_r(f):=\|f\|_{\Hoel^r}$ if $r\in(0,\infty)\setminus\Nat$.

\medskip

\begin{lem*}
For any $f\in C^\infty_b(\bbR^d)$ and $r>0$, we have the estimate
\begin{align}\label{app: Taylor on U alpha}
\left|(\varphi_\lambda f)(y)-\sum_{|k|<r}\frac{\partial^k(\varphi_\lambda f)(x)}{k!}\big(\psi_\lambda(y)-\psi_\lambda(x)\big)^k\right|\lesssim B_r(f)|y-x|^r.
\end{align}
\end{lem*}

\medskip

\begin{Dem}
For $x,y\in U_\lambda$, since $(\psi_\lambda(y)-\psi_\lambda(x))^k=(y-x)^k$, equation \eqref{app: Taylor on U alpha} is just a usual Taylor expansion. For $x,y\notin U_\lambda$, the left hand side of \eqref{app: Taylor on U alpha} is equal to $0$. Let $y\in U_\lambda$ and $x\notin U_\lambda$. Then the left hand side of \eqref{app: Taylor on U alpha} is equal to $|(\varphi_\lambda f)(y)|$. By assumptions, we have
$$
\big|(\varphi_\lambda f)(y)\big|\lesssim\|f\|_{L^\infty} \, d\big(y,U_\lambda^c\big)^r  \le  \|f\|_{L^\infty}\,|y-x|^r,
$$
with an implicit constant in the first inequality depending only on $r$. We have the same estimate when $x\in U_\lambda$ and $y\notin U_\lambda$.
\end{Dem}

\smallskip

Define now $E := \big\{0,\frac14,\frac12,\frac34\big\}^d$ and set
$$
\phi_e := \sum_{\lambda\equiv e\,\text{mod $\Int^d$}}\varphi_\lambda,  \qquad   
x_e^i := \sum_{\lambda\equiv e\,\text{mod $\Int^d$}}\psi_\lambda^i.
$$
Since they are sums of functions with disjoint supports, we have $\phi_e, x_e^i\in C_b^\infty(\bbR^d)$.

\smallskip

\noindent \textbf{\textsf{Proof of Proposition \ref{app: prop PTE}} -- } For $x=(x_i)_{i=1}^d$, we define $|x|_\infty:=\sup_{i=1,\dots,d}|x_i|$. If $|y-x|_\infty\ge\frac12$, then the left hand side of \eqref{app: patched Taylor expansion} is bounded by
$CB_r(f)\lesssim CB_r(f)|y-x|_\infty$, with some constant $C$ depending on $\Phi_{\bar\lambda}$ and $\psi_{\bar\lambda}^i$. If $|y-x|_\infty<\frac12$, then there is no pair $(\lambda,\lambda')$ such that $\lambda\neq\lambda'$ and $(x,y)\in V_\lambda\times V_{\lambda'}$. Hence there exists $\lambda\in\Lambda$ such that the left hand side of \eqref{app: patched Taylor expansion} is equal to that of \eqref{app: Taylor on U alpha}, so the required estimate follows.   \hfill $\rhd$

\bigskip

\noindent \textbf{\textsf{Proof of Proposition \ref{PropPolynomialLift} -- }} We need to show that the component of $\bsf(y)-\widehat{{\sf g}_{yx}}\bsf(x)$ on ${\bfX}_e^k$ is no greater than a constant multiple of $B_r(f) \, |y-x|^{r-|k|}$. Note that $\widehat{{\sf g}_{yx}}\bsf(x)$ is given by
$$
\widehat{{\sf g}_{yx}}\bsf(x) = \sum_{e\in E}\sum_{|\ell+m|<r}\frac{\partial^{\ell+m}(\phi_e f)(x)}{\ell!m!}\big(x_e(y) - x_e(x)\big)^\ell \, {\bfX}_e^m.
$$
The $\unit$-coefficient of $\bsf(y)-\widehat{{\sf g}_{yx}}\bsf(x)$ is
$$
\left|\sum_{e\in E}(\phi_e f)(y)-\sum_{e\in E}\sum_{|k|<r}\frac{\partial^k \big(\phi_e f)(x)}{k!}(x_e(y) - x_e(x)\big)^k\right|\lesssim B_r(f) \, |y-x|^r.
$$
This is  estimate \eqref{app: patched Taylor expansion}. For the $\bfX_e^{k}$-coefficient of $\bsf(y)-\widehat{{\sf g}_{yx}}\bsf(x)$, with $k\neq0$, one has
\begin{align}\label{app: Taylor for Phi_lambda^k}
\frac1{k!}\left|\partial^k(\phi_e f)(x)-\sum_{|\ell|<r-|k|}\frac{\partial^\ell\partial^k(\phi_e f)(x)}{\ell!}\big(x_e(y) - x_e(x)\big)^\ell\right| \lesssim B_r(f) \, |y-x|^{r-|k|}.
\end{align}
This is shown by the similar argument to the proof of the estimate \eqref{app: Taylor on U alpha}.   \hfill $\rhd$

\bigskip
\bigskip

\noindent \textcolor{gray}{$\bullet$} {\sf I. Bailleul} -- Univ. Rennes, CNRS, IRMAR - UMR 6625, F-35000 Rennes, France.   \\
\noindent {\it E-mail}: ismael.bailleul@univ-rennes1.fr   

\medskip

\noindent \textcolor{gray}{$\bullet$} {\sf M. Hoshino} --  Faculty of Mathematics, Kyushu University, Japan   \\
{\it E-mail}: hoshino@math.kyushu-u.ac.jp

\end{document}